\documentclass[twoside,11pt,hidelinks]{amsart}

\usepackage[utf8]{inputenc}
\usepackage[english]{babel}
\usepackage[fixlanguage]{babelbib}
\usepackage{amssymb, amsfonts, amsmath, amsthm, amscd, mathtools, thmtools, tikz, xcolor, hyperref, a4wide}
\usepackage[numbers,square]{natbib}
\usepackage[all]{xy}
\usepackage{ xfrac, calc, indentfirst, cleveref, aligned-overset, easyReview, enumitem, float, tikz-cd, microtype, tensor, epsfig, stackengine, scalerel, wasysym, color, graphicx }

\newcounter{contador}
\numberwithin{contador}{section}

\newtheorem{theorem}[contador]{Theorem}
\newtheorem{prop}[contador]{Proposition}
\newtheorem{lemma}[contador]{Lemma}
\newtheorem{corollary}[contador]{Corollary}

\theoremstyle{definition}
\newtheorem{defi}[contador]{Definition}
\newtheorem{obs}[contador]{Remark}
\newtheorem{exe}[contador]{Example}

\newcommand{\Not}[2][0]{	
	\setcounter{enumi}{#1}
	\renewcommand{\theenumi}{#2\arabic{enumi}}
	\renewcommand{\labelenumi}{(\theenumi)}
	\setlength{\itemindent}{\widthof{#2}}
	\setlength{\itemsep}{4pt}
}

\allowdisplaybreaks 

\title[An Ehresmann-Schein-Nambooripad-type theorem]{An Ehresmann-Schein-Nambooripad-type theorem for left restriction semigroupoids}
\author[Haag, Lautenschlaeger and Tamusiunas]{Rafael Haag, Wesley G. Lautenschlaeger and Thaísa Tamusiunas$^*$}
\address{Instituto de Matem\'{a}tica, Universidade Federal do Rio Grande do Sul,  Av. Bento Gon\c{c}alves, 9500, 91509-900. Porto Alegre-RS, Brazil}
\email{rafaelpetasny@gmail.com}
\email{wesleyglautenschlaeger@gmail.com}
\email{thaisa.tamusiunas@gmail.com}
\thanks{$^*$ Corresponding author}
\date{} 

\begin{document}

    
    \begin{abstract}
       We introduce the concept of locally inductive constellations and establish isomorphisms between the categories of left restriction semigroupoids and locally inductive constellations. This construction offers an alternative to the celebrated Ehresmann–Schein–Nambooripad (ESN) Theorem and, in particular, generalizes results for one-sided restriction semigroups. We also obtain ESN-type theorems for one-sided restriction categories and inverse semigroupoids. 
    \end{abstract}

    \maketitle

\section{Introduction}
    The Ehresmann-Schein-Nambooripad Theorem (ESN Theorem for short), named after Charles Ehresmann, Boris Schein and Komana Nambooripad for their contribution to the theory of semigroups and partial functions, states that:
    
    \begin{theorem} \cite[Theorem 4.1.8]{lawson1998inverse} \label{ESN-Theorem}
        The category of inverse semigroups and premorphisms is isomorphic to the category of inductive groupoids and ordered functors; and the category of inverse semigroups and morphisms is isomorphic to the category of inductive groupoids and inductive functors.
    \end{theorem}
    
    Prior to its appearance in Lawson's book, the ESN Theorem has been generalized to regular semigroups by Nambooripad \cite{nambooripad1974} and to concordant semigroups by Sheena Armstrong \cite{armstrong1988}, and more recently to the class of weakly $U$-regular semigroups by Yanhui Wang \cite{wang2016}. In 1991, the ESN Theorem was generalized by Lawson in another direction, concluding the following:

    \begin{theorem} \cite[Theorem 4.24]{lawson1991}
        The category of Ehresmann semigroups and admissible morphisms is isomorphic to the category of Ehresmann categories and strongly ordered functors.
    \end{theorem}

    Lawson's approach to the ESN Theorem for Ehresmann semigroups was generalized to weakly B-orthodox semigroups by Y. Wang \cite{wang2014}, to locally restriction $P$-restriction semigroups by Shoufeng Wang \cite{wang2019}, to locally Ehresmann $P$-Ehresmann semigroups by S. Wang \cite{wang2020}, to localizable semigroups by FitzGerald and Kinyon \cite{fitz2021}, to DRC-semigroups by S. Wang \cite{wang2022}, and to cat-semigroups satisfying the left match-up condition by Tim Stokes \cite{stokes2023}. The relation between some of those classes of semigroups was studied by S. Wang in \cite{wang2023d}. On the other hand, Darien DeWolf and Dorette Pronk generalized Theorem \ref{ESN-Theorem} in the following way:

    \begin{theorem} \cite[Corollary 3.18 and Theorem 3.16]{dewolf2018ehresmann} \label{multi-ESN}
        The category of inverse semicategories and semifunctors is isomorphic to the category of locally inductive groupoids and locally inductive functors; and the category of inverse categories and functors is isomorphic to the category of locally complete inductive groupoids and locally inductive functors.
    \end{theorem}

    Semicategories are also known as Tilson semigroupoids, graphed semigroupoids or categorical semigroupoids, and constitute a proper subclass of semigroupoids introduced by Exel in \cite{exel2011semigroupoid}.

    The previously mentioned generalizations of Theorem \ref{ESN-Theorem} are considered two-sided ESN-type theorems, since they consider structures together with two unary operations and/or a set of idempotent elements with two partial orders. For instance, an inverse semigroup can be regarded as a two-sided restriction semigroup by defining $s^+ = ss^{-1}$ and $s^\ast = s^{-1}s$, and in this case the partial orders
        $$ s \leq_l t \iff \exists e \in E(\mathcal{S}) \colon s = et \quad\text{and}\quad s \leq_r t \iff \exists e \in E(\mathcal{S}) \colon s = te, $$
    coincide. In order to develop a one-sided ESN-type theorem, Victoria Gould and Christopher Hollings introduced the notion of left constellations and radiants, a one-sided generalization of categories and functors. In \cite[Theorem 4.13]{gould2009restriction} and \cite[Theorem 6.6]{gould2011actions}, Gould and Hollings obtained the following:

    \begin{theorem} \label{left-ESN}
        The category of left restriction semigroups and (2,1)-morphisms is isomorphic to the category of inductive left constellations and ordered radiants; and the category of left restriction semigroups and strong premorphisms is isomorphic to the category of inductive left constellations and ordered preradiants.
    \end{theorem}

    Gould and Hollings one-sided ESN-type theorem was generalized to left congruence D-semigroups by Stokes \cite{stokes2017d}. The relation between several classes of one-sided generalizations of inverse semigroups is also considered in \cite{wang2023d}. We emphasize that Theorem \ref{left-ESN} is not a generalization of Theorem \ref{ESN-Theorem}, but an alternative to it. In fact, \cite[Corollary 5.4]{gould2009restriction} relates the category of inverse semigroups and morphisms to the category of inductive left constellations with right inverses and ordered radiants, which is isomorphic (but not equal) to the category of inductive groupoids and inductive functors. A relation between constellations and categories is considered in \cite{gould2017constcat}.

    The aim of this paper is to establish an ESN-type theorem for the class of left restriction semigroupoids. This class includes the left restriction semigroups studied in \cite{gould2009restriction}, the left restriction categories introduced by Cockett and Lack in \cite{cockett2002restriction}, and the inverse semigroupoids considered in \cite{dewolf2018ehresmann}. In particular, we extend Theorem \ref{left-ESN}, provide a version for left restriction categories, which has not previously appeared in the literature, and introduce an alternative formulation, which we refer to as the “one-sided ESN-type theorem” for inverse semigroupoids.

   The paper is organized as follows. Section 2 recalls the definition and fundamental properties of left restriction semigroupoids and their morphisms. Section 3 introduces the notions of locally inductive constellations, their morphisms, and their properties. In Section 4 we present the first part of the ESN-type theorem by establishing an isomorphism between the category of left restriction semigroupoids with morphisms and the category of locally inductive constellations with inductive radiants. Section 5 introduces the Szendrei expansion of a locally inductive constellation and establishes its universal property. In Section 6 we prove the second part of the ESN-type theorem, showing an isomorphism between the category of left restriction semigroupoids with premorphisms and the category of locally inductive constellations with inductive preradiants. Finally, Section 7 applies our main result to the particular cases of left restriction categories, left restriction semigroups, and inverse semigroupoids regarded as left restriction structures.

\section{Left Restriction Semigroupoids}
 In this section we recall the definition of a left restriction semigroupoid, as introduced in \cite{rsgpdexpansion}. We also introduce the category of left restriction semigroupoids and morphisms, as well as the category of left restriction semigroupoids and premorphisms.

    \begin{defi} \label{def:semigroupoid} \cite[Definition 2.1]{exel2011semigroupoid}
        A \textit{semigroupoid} is a triple $S = (S,S^{(2)},\star)$ such that $S$ is a set, $S^{(2)}$ is a subset of $S \times S$, and $\star \colon S^{(2)} \to S$ is an operation which is associative in the following sense: if $r,s,t \in S$ are such that either
        \begin{enumerate} \Not{s}
            \item $(s,t) \in S^{(2)}$ and $(t,r) \in S^{(2)}$, or \label{s1}
            \item $(s,t) \in S^{(2)}$ and $(s \star t,r) \in S^{(2)}$, or \label{s2}
            \item $(t,r) \in S^{(2)}$ and $(s,t \star r) \in S^{(2)}$, \label{s3}
        \end{enumerate}
        then all of $(s,t)$, $(t,r)$, $(s \star t,r)$ and $(s,t \star r)$ lie in $S^{(2)}$, and $(s \star t) \star r = s \star (t \star r)$.
    \end{defi}

    At times, we will refer to $(s,t) \in S^{(2)}$ as “$st$ is defined", omitting the composition symbol and using concatenation to represent it. For each $s \in S$ consider the sets
        $$ {^sS} = \{ t \in S \colon (t,s) \in S^{(2)} \} \quad\text{and}\quad S^s = \{t \in S \colon (s,t) \in S^{(2)}\}. $$
    Then the associativity of $S$ implies that $S^{st} = S^t$ and ${^{st}S} = {^sS}$, whenever $st$ is defined.

    \begin{defi} \label{def:lr-sg} \cite[Definition 2.5]{rsgpdexpansion}
        A \textit{left restriction semigroupoid} is a pair $(S,+)$ in which $S$ is a semigroupoid and $s \mapsto s^*$ is a unary operation such that
        \begin{enumerate}\Not{lr}
            \item for all $s \in S$, $s^+s$ is defined and $s^+s = s$; \label{lr1}
            \item $s^+t^+$ is defined if and only if $t^+s^+$ is defined, and in this case $s^+ t^+ = t^+s^+$; \label{lr2}
            \item if $s^+t$ is defined, then $(s^+t)^+ = s^+t^+$; \label{lr3}
            \item if $st$ is defined, then $st^+ = (st)^+s$. \label{lr4}
        \end{enumerate}
    In this case, we denote $S^+ = \{ s^+ \colon s \in S \}$. We will usually write $S$ instead of the pair $(S,+)$ when no confusion can arise.
    \end{defi}

    By \cite[Remark 2.6]{rsgpdexpansion}, the left restriction semigroupoid axioms are well defined. Furthermore, it follows from \cite[Proposition 2.15]{rsgpdexpansion} that every left restriction semigroupoid is a categorical semigroupoid. From \eqref{lr1}, we have that ${^tS} = {^{t^+t}S} = {^{t^+}S}$, for every $t \in S$. That is, $st$ is defined if and only if $st^+$ is defined. We recall that a left restriction semigroupoid has a natural partial order, defined by
        $$ s \leq t \iff s^+t \text{ is defined and } s^+t = s. $$
    Or equivalently,
        $$ s \leq t \iff \exists e \in S^+ \colon et \text{ is defined and } et = s. $$
    For a more detailed discussion on left restriction semigroupoids, see \cite{rsgpdexpansion}.

    \begin{lemma} \label{lema:lr} \cite[Lemmas 2.12 and 2.14]{rsgpdexpansion}
        Let $S$ be a left restriction semigroupoid. Then the following statements hold.
        \begin{itemize}
            \item[(a)] For every $e \in S^+$, $ee$ is defined and $ee = e$.
            \item[(b)] For every $e \in S^+$, we have $e^+ = e$.
            \item[(c)] If $st$ is defined, then $(st^+)^+ = (st)^+ = (st)^+ s^+$.
        \end{itemize}
    \end{lemma}

   Recall that a category is a quintuple $(C_0,C_1,D,R,\circ)$, where $C_0$ is a family of \textit{objects}, $C_1$ is a family of \textit{morphisms}, $D,R \colon C_1 \to C_0$ are the \textit{domain} and \textit{range} functions, and $\circ$ is a partially defined, unital, and associative operation on $C_1$. The set of composable pairs is given by
        $$ C_1^{(2)} = \{ (f,g) \in C_1 \times C_1 \colon D(f) = R(g) \}. $$
Moreover, for any $(f, g) \in C_1^{(2)}$, the composition satisfies $D(fg) = D(g)$ and $R(fg) = R(f)$.

We now proceed to define the categories of left restriction semigroupoids.

    \begin{defi}
        A function $\varphi \colon S \to T$ between left restriction semigroupoids is called a \textit{morphism} if the following conditions are satisfied:
        \begin{enumerate} \Not{rm}
            \item if $st$ is defined, then $\varphi(s)\varphi(t)$ is defined and $\varphi(st) = \varphi(s)\varphi(t)$; \label{rm1}
            \item $\varphi(s^+) = \varphi(s)^+$, for every $s \in S$. \label{rm2}
        \end{enumerate}

        In other words, a morphism between left restriction semigroupoids preserves both the binary and the unary operations.
    \end{defi}

    \begin{prop} \label{prop:sgpd-cat-1}
        Left restriction semigroupoids (as objects) and morphisms form a category, where the domain and range are given by the domain and range of functions, and composition is given by the usual composition of functions.

        \begin{proof}
            We only need to verify that the identity map $\mathrm{id}_S \colon S \to S$ is a morphism, and that the composition of two morphisms is again a morphism. It is straightforward to check that $\mathrm{id}_S$ satisfies the axioms of a morphism.

        To verify that the composition is a morphism, let $\varphi \colon S \to T$ and $\psi \colon T \to L$ be morphisms. Take elements $s, t \in S$ such that $st$ is defined. Since $\varphi$ is a morphism, condition \eqref{rm1} ensures that $\varphi(s)\varphi(t)$ is defined. Similarly, by applying \eqref{rm1} to $\psi$, we conclude that $\psi(\varphi(s))\psi(\varphi(t))$ is also defined. Hence,
                $$ \psi(\varphi(s))\psi(\varphi(t)) = \psi(\varphi(s)\varphi(t)) = \psi(\varphi(st)). $$
            Furthermore, applying \eqref{rm2} to $\varphi$ and to $\psi$, we obtain:
                $$ \psi(\varphi(s^+)) = \psi(\varphi(s)^+) = \psi(\varphi(s))^+. $$
            Therefore, the composition $\psi \circ \varphi$ is a morphism.
        \end{proof}
    \end{prop}

    \begin{defi}
        A function $\varphi \colon S \to T$ between left restriction semigroupoids is called a \textit{premorphism} if the following conditions are satisfied:
        \begin{enumerate} \Not{pm}
            \item If $st$ is defined, then $\varphi(s)\varphi(t)$ and $\varphi(s)^+ \varphi(st)$ are defined and $\varphi(s)\varphi(t) = \varphi(s)^+ \varphi(st)$.\label{pm1}
            \item $\varphi(s)^+ \leq \varphi(s^+)$, for every $s \in S$. That is, $\varphi(s)^+ = \varphi(s)^+\varphi(s^+)^+$. \label{pm2}
        \end{enumerate}
    \end{defi}

    In some contexts, a function $\varphi$ is called a premorphism if it satisfies $\varphi(s)\varphi(t) \leq \varphi(st)$ together with condition \eqref{pm2}. However, \eqref{pm1} is stronger than the inequality $\varphi(s)\varphi(t) \leq \varphi(st)$. Accordingly, a function satisfying both \eqref{pm1} and \eqref{pm2} is also referred to as a \textit{strong} premorphism, as in \cite{gould2011actions}. In the present paper, we are concerned only with this stronger notion. Therefore, we will simply refer to strong premorphisms as premorphisms.

    \begin{lemma} \label{lema:pm} \cite[Lemma 4.4]{rsgpdexpansion}
        Let $\varphi \colon S \to T$ be a premorphism. Then the following statements hold.
        \begin{itemize}
            \item[(a)] $\varphi(e) \in T^+$, for every $e \in S^+$. \label{lema:pm-a}
            
            \item[(b)] If $s,t \in S$ and $s \leq t$, then $\varphi(s) \leq \varphi(t)$. \label{lema:pm-b}
        \end{itemize}
    \end{lemma}

    \begin{prop} \label{prop:sgpd-cat-2}
        Left restriction semigroupoids (as objects) and premorphisms (as morphisms) form a category, where the domain and range are given by the domain and range of functions, and $\circ$ is the usual composition of functions.

        \begin{proof}
            It suffices to verify that $id_S$ is a premorphism, and that the composition $\psi \circ \varphi \colon S \to L$ is a premorphism whenever $\varphi \colon S \to T$ and $\psi \colon T \to L$ are premorphisms. It is straightforward that $id_S$ satisfies the conditions of a premorphism.

            For the composition, let $s,t \in S$ be such that $st$ is defined. Since $\varphi$ is a premorphism, it follows from condition \eqref{pm1} that $\varphi(s)\varphi(t)$ and $\varphi(s)^+ \varphi(st)$ are defined. Moreover, since $\psi$ is a premorphism and $\varphi(s)\varphi(t)$ is defined, condition \eqref{pm1} ensures that $\psi(\varphi(s))\psi(\varphi(t))$ and $\psi(\varphi(s))^+ \psi(\varphi(s)\varphi(t))$ are also defined. In this case, we have:
            \begin{align*}
                \psi(\varphi(s))\psi(\varphi(t)) &= \psi(\varphi(s))^+ \psi(\varphi(s)\varphi(t)) & \eqref{pm1}, \psi \\
                &= \psi(\varphi(s))^+ \psi(\varphi(s)^+\varphi(st)) & \eqref{pm1}, \varphi \\
                &= \psi(\varphi(s))^+ \psi(\varphi(s)^+)^+ \psi(\varphi(s)^+\varphi(st)) & \eqref{pm2}, \psi \\
                &= \psi(\varphi(s))^+ \psi(\varphi(s)^+) \psi(\varphi(st)) & \eqref{pm1}, \psi \\
                &= \psi(\varphi(s))^+ \psi(\varphi(s)^+)^+ \psi(\varphi(st)) & \ref{lema:pm}(a), \ref{lema:lr}(b) \\
                &= \psi(\varphi(s))^+ \psi(\varphi(st)). & \eqref{pm2}, \psi
            \end{align*}
            
            On the other hand, for every $s \in S$, we have $\varphi(s^+) \leq \varphi(s)^+$ by \eqref{pm2} for $\varphi$. Therefore,
            \begin{align*}
                \psi(\varphi(s^+)) \leq \psi(\varphi(s)^+) \leq \psi(\varphi(s))^+,
            \end{align*} where the first inequality follows from Lemma \ref{lema:pm}(b), and the second from condition \eqref{pm2} for $\psi$. Hence, $\psi \circ \varphi$ satisfies \eqref{pm2}, and thus is a premorphism.
        \end{proof}
    \end{prop}

    \begin{obs}
        Every restriction morphism is a premorphism. In fact, if $\varphi$ is a restriction morphism and $st$ is defined, then
            $$ \varphi(s)\varphi(t) \overset{\eqref{lr1}}{=} \varphi(s)^+ \varphi(s) \varphi(t) \overset{\eqref{rm1}}{=} \varphi(s)^+ \varphi(st). $$
       Furthermore, for every $s \in S$, we have $\varphi(s^+) = \varphi(s)^+$, which is a special case of the inequality $\varphi(s^+) \leq \varphi(s)^+$. Hence, the category of left restriction semigroupoids with restriction morphisms is a subcategory of the one with premorphisms.
    \end{obs}

\section{Locally Inductive Left Constellations} In this section we introduce locally inductive left constellations (li-constellations) and discuss some of their properties.

    Let $X$ be a set endowed with a partially defined binary operation $\star$. An element $e \in X$ is said to be an \textit{idempotent} if $e \star e$ is defined and $e \star e = e$. We denote the set of idempotent elements of $X$ by $E(X)$.

    \begin{defi} \label{def:constellation} \cite[Definition 2.1]{gould2009restriction}
        A \textit{left constellation} is a quadruple $T = (T,T^{(2)},\star,+)$, where $T$ is a set, $T^{(2)} \subseteq T \times T$, $\star \colon T^{(2)} \to T$ is a binary operation, and $+ \colon T \to T$ is a unary operation, satisfying the following axioms:
        \begin{enumerate} \Not{c}
            \item $(x,y),(y,z) \in T^{(2)}$ if and only if $(y,z),(x,y \star z) \in T^{(2)}$; \label{c1}
            \item if $(x,y),(y,z) \in T^{(2)}$, then $(x \star y,z) \in T^{(2)}$, and in this case $x \star (y \star z) = (x \star y) \star z$; \label{c2}
            \item if $e \in T^+ = \{ x^+ \colon x \in T \}$, then $(e,x) \in T^{(2)}$ and $e \star x = x$ if and only if $e=x^+$; \label{c3}
            \item if $e \in T^+$ and $(x,e) \in T^{(2)}$, then $x \star e = x$. \label{c4}
        \end{enumerate}
    \end{defi}

    Analogously to semigroupoids, we may refer to $(x,y) \in T^{(2)}$ as ``$xy$ is defined", omitting the composition symbol and using concatenation to represent it. For each $x \in T$, we denote
        $$ {^xT} = \{ y \in T \colon (y,x) \in T^{(2)}\} \quad\text{and}\quad T^x = \{ y \in T \colon (x,y) \in T^{(2)} \}. $$
    In this case, if $xy$ is defined, then condition \eqref{c1} implies that ${^xT} = {^{xy}T}$, and \eqref{c2} implies that $T^y \subseteq T^{xy}$. For a more detailed discussion on left constellations, see \cite{gould2009restriction}. The next result incorporates \cite[Lemma 2.2]{gould2009restriction} and \cite[Lemma 2.3]{gould2009restriction}.

    \begin{lemma} \label{lema:constellation}
        Let $T$ be a left constellation. The following statements hold.
        \begin{itemize}
            \item[(a)] If $e \in T^+$, then $ee$ is defined, $ee = e$ and $e^+ = e$. That is, $T^+ \subseteq E(T)$.
            \item[(b)] If $xy$ is defined, then $(xy)^+ = x^+$.
            \item[(c)] $xy$ is defined if and only if $xy^+$ is defined.
        \end{itemize}
    \end{lemma}

    Before introducing li-constellations, we recall the notion of inductive constellations as presented in \cite{gould2009restriction}, along with some definitions concerning partially ordered sets.

    Let $(X,\leq)$ be a partially ordered set. There is a natural equivalence relation on $X$, defined by $x \sim y$ if and only if there exist elements $z_1, \dots, z_n \in X$ such that
        $$ x = z_1, \quad z_n = y \quad\text{and}\quad z_i \leq z_{i+1} \text{ or } z_{i+1} \leq z_i, \ \forall i=1,\dots,n-1. $$
    Since $\sim$ is an equivalence relation, the set $X$ can be uniquely expressed as the disjoint union of its $\sim$-equivalence classes, which we refer to as the \textit{connected components} of $X$. Given $x \in X$, we denote its connected component by
        $$ \omega(x) = \{ y \in X \colon y \sim x \}. $$
    Each connected component of $X$ is itself a partially ordered set with the same order inherited from $X$.

    \begin{defi}
        Let $(X,\leq)$ be a partially ordered set. We say that $(X,\leq)$ is:
        \begin{itemize}
            \item[(1)] a \textit{meet-semilattice} if, for every pair $x,y \in X$, there exists $z \in X$ satisfying:
            \begin{align*}
                [z \leq x,y] \quad\text{and}\quad [z' \leq x, y \implies z' \leq z].
            \end{align*} In this case, we write $z = x \wedge y$ and call $z$ the \emph{meet} of $x$ and $y$.

            \item[(2)] a \textit{local meet-semilattice} if each connected component $\omega(x)$ of $X$ is a meet-semilattice.
        \end{itemize}
    \end{defi}

    Since the partial order $\leq$ is antisymmetric, the element $x \wedge y$, whenever it exists, is unique. Thus, $(X, \leq)$ is a meet-semilattice if and only if the operation $\wedge \colon X \times X \to X$ endows $X$ with a semigroup structure in which $x \wedge y = y \wedge x$ and $x \wedge x = x$ for all $x, y \in X$. In other words, meet-semilattices are precisely the commutative semigroups in which every element is idempotent. Local meet-semilattices are disjoint unions of such semigroups, and therefore are semigroupoids.

    \begin{defi} \cite[Definition 3.3]{gould2009restriction} \label{defi:i-constellation}
        Let $T = (T,T^{(2)},+,\star)$ be a left constellations and $\leq$ be a partial order on $T$. We say that $(T,\leq)$ is \textit{inductive} if, for $x,x',y,y' \in T$, the following conditions are satisfied:
        \begin{enumerate} \Not{O}
            \item if $x \leq y$, $x' \leq y'$ and $xx'$ and $yy'$ are defined, then $xx' \leq yy'$; \label{gould-const-1}
            \item if $x \leq y$, then $x^+ \leq y^+$; \label{gould-const-2}
            \item if $e \in T^+$ and $e \leq x^+$, then there is a unique element $e|x \in T$ satisfying $e|x \leq x$ and $(e|x)^+ = e$, which we call the \textit{restriction} of $x$ to $e$; \label{gould-const-3}
            \item for every $e \in T^+$, there is a maximum element $x|e \in T$ such that $x|e \leq x$ and $(x|e)e$ is defined, which we call the \textit{corestriction} of $x$ to $e$; \label{gould-const-4}
            \item if $e \in T^+$ and $xy$ is defined, then $((xy)|e)^+ = (x|(y|e)^+)^+$; \label{gould-const-5}
            \item if $e,f \in T^+$ and $f \leq e$, then the restriction $f|e$ coincides with the corestriction $f|e$; \label{gould-const-6}
            \item $(T^+,\leq)$ is a meet-semilattice in which the meet operation is given by corestriction. \label{gould-const-7}
        \end{enumerate}
    \end{defi}

    The first six axioms define an \emph{ordered left constellation}, while the six axioms together with the last one characterizes an \emph{inductive left constellation}, which is a particular case of an ordered one. Inductive left constellations were introduced in \cite{gould2009restriction} to establish isomorphisms with categories of left restriction semigroups (see \cite[Theorem 4.13]{gould2009restriction} and \cite[Theorem 6.6]{gould2011actions}). To obtain isomorphisms with the categories of left restriction semigroupoids described in Propositions \ref{prop:sgpd-cat-1} and \ref{prop:sgpd-cat-2}, a broader class of partially ordered left constellations is required.
    
    \begin{defi}
        Let $T = (T,T^{(2)},+,\star)$ be a left constellation and $\leq$ be a partial order in $T$. We say that $(T,\leq)$ is \textit{locally inductive} if, for $x,x',y,y' \in T$, the following conditions hold:
        \begin{enumerate} \Not{wo}
            \item if $x \leq y$, $x' \leq y'$ and $xx'$ and $yy'$ are defined, then $xx' \leq yy'$; \label{wo1}
            
            \item if $x \leq y$, then $x^+ \leq y^+$; \label{wo2}
            
            \item if $e \in T^+$ and $e \leq x^+$, then there is a unique element $e|x \in T$ such that $e|x \leq x$ and $(e|x)^+ = e$, which we call the \textit{restriction} of $x$ to $e$; \label{wo3}
            
            \item for every $e \in T^+$, the set $\{ y \in T \colon y \leq x \text{ and $ye$ is defined} \}$ is either empty or has a maximum element. In the latter case, we denote this maximum element by $x|e$, refer to it as the \textit{corestriction} of $x$ to $e$, and write $x|e \neq \emptyset$ to indicate that it exists; \label{wo4}

            \item if $e \in T^+$ and $xy$ is defined, then $(xy)|e \neq \emptyset$ if and only if $y|e \neq \emptyset$; \label{wo5}

            \item if $e,f \in T^+$ and $f \leq e$, then $x|e \neq \emptyset$ if and only if $x|f \neq \emptyset$; \label{wo6}

            \item if $e \in T^+$, $xy$ is defined and $(xy)|e \neq \emptyset$, then $((xy)|e)^+ = (x|(y|e)^+)^+$; \label{wo7}

            \item if $e,f \in T^+$ and $e \leq f$, then the restriction $e|f$ and the corestriction $e|f$ coincide;\label{wo8}

            \item $(T^+,\leq)$ is a local meet-semilattice with meet operation given by corestriction. \label{wo9}
        \end{enumerate}
    \end{defi}

    \begin{obs} \label{obs:li-constellations}
        Conditions \eqref{wo7} and \eqref{wo8} are well defined.

Indeed, let $e \in T^+$ and $x, y \in T$ be such that $xy$ is defined, and suppose that $(xy)|e \neq \emptyset$. By \eqref{wo5}, we have $y|e \neq \emptyset$, and since $y|e \leq y$, it follows from \eqref{wo2} that $(y|e)^+ \leq y^+$. On the other hand, since $xy$ is defined, Lemma~\ref{lema:constellation}(c) implies that $xy^+$ is defined. Therefore,
            $$ x \in \{ z \in T \colon z \leq x \text{ and $zy^+$ is defined} \}, $$
        which shows that $x|y^+ \neq \emptyset$. Since $(y|e)^+ \leq y^+$, it follows from \eqref{wo6} that $x|(y|e)^+ \neq \emptyset$. Hence, condition \eqref{wo7} is well defined. 
        
        Now, suppose $e, f \in T^+$ are such that $e \leq f$. Then the restriction $e|f$ exists by definition. Moreover, Lemma~\ref{lema:constellation}(a) ensures that
            $$ e \in \{ z \in T \colon z \leq e \text{ and $ze$ is defined} \}. $$
       which implies $e|e \neq \emptyset$. Since $e \leq f$, it follows from \eqref{wo6} that $e|f \neq \emptyset$. Thus, condition \eqref{wo8} is also well defined.
    \end{obs}

From now on we will refer to locally inductive left constellations simply as \emph{li-constellations}.

    \begin{prop} \label{prop:i-constellations}
        Inductive left constellations are precisely those li-constellations in which the corestriction of $x$ to $e$ exists for every $x \in T$ and $e \in T^+$.

        \begin{proof}
            In general, conditions \eqref{gould-const-1}, \eqref{gould-const-2}, \eqref{gould-const-3} and \eqref{gould-const-6} correspond to conditions \eqref{wo1}, \eqref{wo2}, \eqref{wo3} and \eqref{wo8}, respectively. Suppose that $x|e \neq \emptyset$ for every $x \in T$ and $e \in T^+$. Then condition \eqref{gould-const-4} and \eqref{gould-const-5} are precisely conditions \eqref{wo4} and \eqref{wo7}, while conditions \eqref{wo5} and \eqref{wo6} are trivially satisfied.

            Now, suppose that $(T,\leq)$ is an inductive left constellation. Then every corestriction is defined and, by condition \eqref{gould-const-7}, $(T^+,\leq)$ is a meet-semilattice with meet operation given by corestriction. Therefore, condition \eqref{wo9} holds and $(T,\leq)$ is a li-constellation.
            
            Conversely, suppose that $(T,\leq)$ is a li-constellation in which $x|e \neq \emptyset$ for every $x \in T$ and $e \in T^+$. From \eqref{wo9}, it follows that $(T^+,\leq)$ is a local meet-semilattice with meet operation given by corestriction. Since the corestriction $e|f$ is defined for every $e,f \in T^+$, $(T^+,\leq)$ must have a single connected component. Therefore, $(T^+,\leq)$ is a meet-semilattice and condition \eqref{gould-const-7} is satisfied. This shows that $(T,\leq)$ is an inductive left constellation.
        \end{proof}
    \end{prop}

    From now on, we fix $T = (T, \leq)$ to be a li-constellation. In what follows, we prove some properties involving restrictions and corestrictions in $T$.

    \begin{prop}[Restriction properties] \label{prop:restriction}
        Let $x,y,z \in T$ and $e,f \in T^+$. Then the following hold.
        \begin{itemize}
            \item[(a)] If $x \leq y$, then $x^+|y = x$.
            \item[(b)] If $e \leq f$, then $ef$ is defined and $e = e|f = ef$.
            \item[(c)] If $e \leq f \leq x^+$, then $e|(f|x) = e|x$ and $e|x \leq f|x$.
            \item[(d)] If $e \leq x^+$, then $ex$ is defined and $e|x = ex$.
        \end{itemize}

        \begin{proof}
            (a) If $x \leq y$, then $x^+ \leq y^+$ by \eqref{wo2}. By \eqref{wo3} there is a unique $x^+|y \in T$ satisfying $x^+|y \leq y$ and $(x^+|y)^+ = x^+$. On the other hand, we have $x \leq y$ by hypothesis and $x^+ = x^+$ trivially. From the uniqueness in \eqref{wo3} we obtain $x^+|y = x$.\\

            (b) If $e \leq f$, then $e = e^+|f = e|f$ by (a) and Lemma \ref{lema:constellation}(a). By \eqref{wo8} we obtain that $e|f \neq \emptyset$ and $e|f = e$. In particular $ef = (e|f)f$ is defined. Since $e \leq f$, $f \leq f$ and $ef$ and $ff$ are defined, it follows from \eqref{wo1} and Lemma \ref{lema:constellation}(a) that $ef \leq ff = f$. Furthermore, from Lemma \ref{lema:constellation}(b) we obtain $(ef)^+ = e^+ = e$. But $e|f$ (as a restriction) is the unique element satisfying $e|f \leq f$ and $(e|f)^+ = e$. Therefore, $e = e|f = ef$.\\

            (c) Since $e \leq x^+$ and $f \leq x^+$, the restriction $e|x$ and $f|x$ are defined. Since $e \leq f = (f|x)^+$, the restriction $e|(f|x)$ is defined. By definition, $e|x$ is the unique element satisfying $e|x \leq x$ and $(e|x)^+ = e$. On the other hand, we have
                $$ e|(f|x) \leq f|x \leq x \quad\text{and}\quad (e|(f|x))^+ = e. $$
            From the uniqueness we obtain $e|(f|x) = e|x$. In particular, $e|x = e|(f|x) \leq f|x$.\\

            (d) Suppose that $e \leq x^+$. Then the restriction $e|x$ is defined and, from (b), we obtain that $ex^+$ is defined. By Lemma \ref{lema:constellation}(c), $ex$ is defined. Since $e \leq x^+$, $x \leq x$ and $ex$ and $x^+x$ are defined, it follows from \eqref{wo1} and \eqref{c3} that $ex \leq x^+x = x$. Furthermore, from Lemma \ref{lema:constellation} we obtain that $(ex)^+ = e^+ = e$. But, $e|x$ is the unique element satisfying $e|x \leq x$ and $(e|x)^+ = e$. Therefore, it must be $e|x = ex$.
        \end{proof}
    \end{prop}

By combining parts (a) and (d) of Proposition \ref{prop:restriction}, we obtain an identity that will be useful in several proofs. For clarity, we state it as a corollary. 

    \begin{corollary} \label{coro:restriction}
        If $x \leq y$, then $x^+y$ is defined and $x = x^+y$.

        \begin{proof}
            If $x \leq y$, then $x^+ \leq y^+$ by \eqref{wo2}. By Proposition \ref{prop:restriction}(a), we have $x = x^+|y$, and applying Proposition \ref{prop:restriction}(d) to the pair $(e, x) = (x^+, y)$, we conclude that $x^+y$ is defined and $x^+|y = x^+y$. Hence, $x = x^+y$.
        \end{proof}
    \end{corollary}

    \begin{prop}[Corestriction properties] \label{prop:corestriction}
        Let $x,y \in T$ and $e,f \in T^+$. Then the following hold.
        \begin{itemize}
            \item[(a)] The product $xy$ is defined if and only if $x|y^+ \neq \emptyset$ and $x|y^+ = x$.
            \item[(b)] If $e \leq f$ and $x|e \neq \emptyset$, then $x|e \leq x|f$.
            \item[(c)] If $e \leq f$ and $x|f \neq \emptyset$, then $(x|f)|e \neq \emptyset$ and $(x|f)|e = x|e$.
        \end{itemize}

        \begin{proof}
            (a) If $xy$ is defined, then by Lemma \ref{lema:constellation}(c), $xy^+$ is also defined. Hence, 
                $$ x \in \{ z \in T \colon z \leq x \text{ and $zy^+$ is defined} \}. $$
            That is, $x|y^+ \neq \emptyset$. Since $x|y^+$ is the maximum element such that $x|y^+ \leq x$ and $(x|y^+)y^+$ is defined, we conclude that $x|y^+ = x$. Conversely, suppose that $x|y^+ \neq \emptyset$ and $x|y^+ = x$. Then $xy^+ = (x|y^+)y^+$ is defined. By Lemma \ref{lema:constellation}(c), it follows that $xy$ is defined.\\
        
            (b) If $e \leq f$ and $x|e \neq \emptyset$, then by \eqref{wo6} we obtain $x|f \neq \emptyset$. By definition, $x|e \leq x$ and $(x|e)e$ is defined. Since $e \leq f$, Corollary \ref{coro:restriction} gives that $ef$ is defined.  Then, using \eqref{c2} and \eqref{c4}, we get that $(x|e)f = [(x|e)e]f$ is defined.  Therefore,
                $$ x|e \in \{ z \in T \colon z \leq x \text{ and $zf$ is defined} \}. $$
            Since $x|f$ is the maximum element in this set, it follows that $x|e \leq x|f$.\\
        
            (c) Suppose that $e \leq f$ and $x|f \neq \emptyset$. Since $x|f \leq x|f$ and $(x|f)f$ is defined, it follows that $(x|f)|f \neq \emptyset$. By \eqref{wo6}, we obtain $(x|f)|e \neq \emptyset$ and $x|e \neq \emptyset$. Observe that, by the definition of corestriction, we have
                $$ (x|f)|e \leq x|f \leq x \quad\text{and}\quad ((x|f)|e)e \text{ is defined}. $$
            Hence,  $(x|f)|e \leq x|e$ since $x|e$ is the maximum element with those properties. On the other hand, since $e \leq f$ and $x|e \neq \emptyset$, it follows from (b) that
                $$ x|e \leq x|f \quad\text{and}\quad (x|e)e \text{ is defined}. $$
            Therefore, $x|e \leq (x|f)|e$ as before. This shows that $x|e = (x|f)|e$.
        \end{proof}
    \end{prop}

    We observe that, up to this point, only conditions \eqref{wo1}–\eqref{wo6} and \eqref{wo8} have been used. In the next result we establish a relation between restriction and corestriction which, together with \eqref{wo7}, allows us to express the corestriction of a product as a product of corestrictions.

    \begin{lemma} \label{lema:corestriction}
        Let $x,y \in T$ and $e \in T^+$.
        \begin{itemize}
            \item[(a)] If $x|e \neq \emptyset$, then $x|e = (x|e)^+|x = (x|e)^+ x$.
            \item[(b)] If $xy$ is defined and $(xy)|e \neq \emptyset$, then $(xy)|e = (x|(y|e)^+)(y|e)$.
        \end{itemize}

        \begin{proof}
            (a) If $x|e \neq \emptyset$, then $x|e \leq x$. Since $(x|e)^+ \leq x^+$ by \eqref{wo2}, the restriction $(x|e)^+|x$ is defined. By Corollary \ref{coro:restriction}, we have $(x|e)^+|x = (x|e)^+ x$. On the other hand, $x|e \leq x$ and $(x|e)^+ = (x|e)^+$, and $(x|e)^+|x$ is the unique element satisfying these conditions. Therefore, $x|e = (x|e)^+|x = (x|e)^+ x$.\\
            
           Suppose that $xy$ is defined and $(xy)|e \neq \emptyset$. Then $x|(y|e)^+$ is defined by condition \eqref{wo7}, and $(x|(y|e)^+)(y|e)^+$ is defined by the definition of corestriction. Hence, $(x|(y|e)^+)(y|e)$ is defined by Lemma \ref{lema:constellation}(c). The result now follows from the calculation:
           
            \begin{align*}
                (x|(y|e)^+)(y|e) &= [(x|(y|e)^+)^+ x][(y|e)^+y] & (a) \\
                &= [[(x|(y|e)^+)^+ x](y|e)^+]y & \eqref{c2} \\
                &= [(x|(y|e)^+)^+ x]y & \eqref{c4} \\
                &= [((xy)|e)^+ x]y & \eqref{wo7} \\
                &= ((xy)|e)^+ [xy] & \eqref{c2} \\
                &= (xy)|e. & (a)
            \end{align*}
        \end{proof}
    \end{lemma}

    The notions of inductive radiants and inductive preradiants introduced below generalize the concepts of ordered radiants \cite[Definition 3.8]{gould2009restriction} and ordered preradiants \cite[Definition 3.6]{gould2011actions}, respectively.

    \begin{defi}
        A function $\varphi \colon T \to L$ between li-constellations is called an \textit{inductive radiant} if, for $x,y \in T$ and $e \in T^+$, the following conditions are satisfied:
        \begin{enumerate} \Not{ir}
            \item if $xy$ is defined, then $\varphi(x)\varphi(y)$ is defined and $\varphi(xy) = \varphi(x)\varphi(y)$; \label{ir1}
            \item $\varphi(x)^+ = \varphi(x^+)$. \label{ir2}
            \item if $x \leq y$, then $\varphi(x) \leq \varphi(y)$; \label{ir3}
            \item if $x|e \neq \emptyset$, then $\varphi(x)|\varphi(e) \neq \emptyset$ and $\varphi(x|e) = \varphi(x)|\varphi(e)$. \label{ir4}
        \end{enumerate}
    \end{defi}

    \begin{defi}
        A function $\varphi \colon T \to L$ between li-constellations is called an \textit{inductive preradiant} if, for $x,y \in T$ and $e \in T^+$, the following conditions are satisfied:
        \begin{enumerate} \Not{ip}
            \item if $xy$ is defined, then $\varphi(x)|\varphi(y)^+ \neq \emptyset$, $\varphi(x)^+|\varphi(xy)^+ \neq \emptyset$ and
                $$ \varphi(x) \otimes \varphi(y) = \varphi(x)^+ \otimes \varphi(xy), $$ \label{ip1}where $s \otimes t := (s|t^+)t$, for every $s,t \in L$;
            \item $\varphi(x)^+ \leq \varphi(x^+)$. \label{ip2}
            \item if $x \leq y$, then $\varphi(x) \leq \varphi(y)$; \label{ip3}
            \item if $x|e \neq \emptyset$, then $\varphi(x)|\varphi(e) \neq \emptyset$ and $\varphi(x|e) = \varphi(x)|\varphi(e)$; \label{ip4}
            \item if $e \leq x^+$, then $\varphi(e|x)^+ = \varphi(e)|\varphi(x)^+$ (as corestriction). \label{ip5}
        \end{enumerate}
    \end{defi}

    \begin{obs}
        We need to verify that conditions \eqref{ip4} and \eqref{ip5} are well defined. 
        
        Assume for a moment that  $e \in T^+$ implies $\varphi(e) \in L^+$. Then the corestriction $\varphi(x)|\varphi(e)$ is defined, so \eqref{ip4} is well defined. Furthermore, if $e \leq x^+$, then $\varphi(e) \leq \varphi(x^+)$ by \eqref{ip3}. Since $\varphi(x)^+ \leq \varphi(x^+)$, we conclude that $\varphi(e)$, $\varphi(x^+)$ and $\varphi(x)^+$ belong to the same connected component of $L^+$. By \eqref{wo9}, the corestriction $\varphi(e)|\varphi(x)^+$ is then defined. Therefore, \eqref{ip5} is also well defined.
    \end{obs}

    In the next result we prove the claim that an inductive preradiant satisfies $\varphi(e) \in L^+$ for every $e \in T^+$, noting that the proof does not use conditions \eqref{ip4} and \eqref{ip5}.

    \begin{prop} \label{prop:preradiant}
        Let $\varphi \colon T \to L$ be an inductive preradiant. Then $\varphi(T^+) \subseteq L^+$.

        \begin{proof}
            If $e \in T^+$, then $e^+ = e$ by Lemma \ref{lema:constellation}(a). From \eqref{ip2}, we have $\varphi(e)^+ \leq \varphi(e^+) = \varphi(e)$. By Corollary \ref{coro:restriction}, it follows that $\varphi(e)^+ = (\varphi(e)^+)^+ \varphi(e)$. Since $(\varphi(e)^+)^+ = \varphi(e)^+$ by Lemma \ref{lema:constellation}(a), we conclude:
                $$ \varphi(e)^+ = (\varphi(e)^+)^+\varphi(e) = \varphi(e)^+\varphi(e) = \varphi(e) \in L^+, $$
            where the last identity follows from \eqref{c3}.
        \end{proof}
    \end{prop}

\section{First part of the ESN-type Theorem} In this section we prove that the category of left restriction semigroupoids with morphisms is isomorphic to the category of li-constellations with inductive radiants. This constitutes the first part of our ESN-type Theorem. The second part will be presented in Section 6, where we examine the category of left restriction semigroupoids with premorphisms.\\

    From now on, $S$ denotes a left restriction semigroupoid. We define a li-constellation $(C(S), \leq)$ as follows: let $C(S) = S$, and set
        $$ C(S)^{(2)} = \{ (s,t) \in S \times S \colon \text{$st^+$ is defined and } st^+ = s \}. $$
    Note that for every pair $(s,t) \in C(S)^{(2)}$, the composition $st$ is defined in $S$, since ${^tS} = {^{t^+}S}$ for all $t \in S$. Thus, we define $\bullet \colon C(S)^{(2)} \to C(S)$ by $s \bullet t = st$. Let $+ \colon C(S) \to C(S)$ be the same function $+$ from $S$, and define a relation $\leq$ on $C(S)$ by:
        $$ s \leq t \iff \text{$st$ is defined in $S$ and $s = s^+t$}. $$
    That is, $\leq$ is the natural partial order of $S$. To distinguish the composable pairs in $S$ and $C(S)$, we write “$st$ is defined” when $(s,t) \in S^{(2)}$, and “$s \bullet t$ is defined” when $(s,t) \in C(S)^{(2)}$.

    \begin{lemma} \label{lema:sgpd-to-const}
        The pair $(C(S), \leq)$ is a li-constellation. Moreover, for $e \in C(S)^+$ and $s \in C(S)$, the restriction $e|s$ is $es$ whenever $e \leq s^+$. As for the corestriction, $s|e$ is defined if and only if $se$ is defined, and in that case we have $s|e = se$.

        \begin{proof}
            First we prove that $(C(S),C(S)^{(2)},\bullet,+)$ is a left constellation. Let $s,t,r \in C(S)$. Then:
            \begin{itemize}
                \item[(a)] $s \bullet t$ and $t \bullet r$ are defined if and only if $st^+ = s$ and $tr^+ = t$.
                \item[(b)] $t \bullet r$ and $s \bullet (t \bullet r)$ are defined if and only if $tr^+ = t$ and $s(tr)^+ = s$.
                \item[(c)] $s \bullet t$ and $(s \bullet t) \bullet r$ are defined if and only if $st^+ = s$ and $(st)r^+ = st$.
            \end{itemize}
            
            Suppose that $t \bullet r$ is defined. Then $tr$ is defined and $tr^+ = t$. Applying $+$ to this identity and using Lemma \ref{lema:lr}(c), we obtain $t^+ = (tr^+)^+ = (tr)^+$. Therefore,
                $$ st^+ = s \iff s(tr)^+ = s. $$
            This proves that (a) is equivalent to (b). Furthermore, if (a) holds, then $(st)r^+ = s(tr^+) = st$. Hence, (a) implies (c). Since the composition in $S$ is associative, we obtain \eqref{c1} and \eqref{c2}.\\

            \noindent\eqref{c3} Let $e \in C(S)^+ = S^+$, and suppose that $e \bullet s$ is defined and $e \bullet s = s$. That is, $es$ is defined, $es^+ = e$ and $es = s$. Then, by \eqref{lr3}, we get
                $$ e = es^+ = (es)^+ = s^+. $$
            On the other hand, $s^+ \bullet s$ is defined since $s^+s^+ = s^+$ by Lemma \ref{lema:lr}(a), and $s^+ \bullet s = s^+s = s$ by \eqref{lr1}. Therefore, $s^+ \in C(S)^+$ is the unique element satisfying \eqref{c3}.\\
    
            \noindent\eqref{c4} Let $e \in C(S)^+$, and suppose that $s \bullet e$ is defined. Then $s \bullet e = se = se^+ = s$ by definition.\\

            This proves that $C(S) = (C(S),C(S)^{(2)},\bullet,+)$ is a left constellation. We now proceed to show that $(C(S),\leq)$ is locally inductive.\\ 
            
            \noindent\eqref{wo1} Suppose that $s \leq t$ and $s' \leq t'$, and that both $s \bullet s'$ and $t \bullet t'$ are defined. Then we have $s = s^+t$ and $s' = s'^+t'$. Since $s^+t$ is defined and $(ss')^+ = (ss')^+ s^+$ by Lemma \ref{lema:lr}(c), it follows that $(ss')^+ tt'$ is also defined. Moreover,
            \begin{align*}
                (ss')^+(tt') &= (ss')^+s^+tt' & \ref{lema:lr}(c) \\
                &= (ss')^+ st' & (s \leq t) \\
                &= ss'^+ t' & \eqref{lr4} \\
                &= ss'. & (s' \leq t')
            \end{align*}
            Therefore, $ss' \leq tt'$.\\

            \noindent\eqref{wo2} If $s \leq t$, then $s = s^+t$. By \eqref{lr3}, it follows that $s^+ = (s^+ t)^+ = s^+t^+$, which implies that $s^+ \leq t^+$.\\

            \noindent\eqref{wo3} Let $e \in C(S)^+$ and $s \in C(S)$ be such that $e \leq s^+$. Then $e = es^+$, and therefore $es$ is defined. Observe that $es$ satisfies:
            \begin{align*}
                (es)^+ s &= (e^+ s)^+ s & \eqref{lema:lr}(b) \\
                &= e^+s^+ s & \eqref{lr3} \\
                &= e^+s & \eqref{lr1} \\
                &= es, & \eqref{lema:lr}(b)
            \end{align*}
            and
            \begin{align*}
                (es)^+ &= (es^+)^+ & \eqref{lema:lr}(c) \\
                &= e^+ & (e \leq s^+) \\
                &= e. & \eqref{lema:lr}(b)
            \end{align*}
            Hence, $es \leq s$ and $(es)^+ = e$. Suppose that $t \in C(S)$ is another element such that $t \leq s$ and $t^+ = e$. Then $t = t^+ s = es$, and therefore $e|s = es$ is the unique element satisfying \eqref{wo3}.\\

            \noindent\eqref{wo4} Let $e \in C(S)^+$ and $s \in C(S)$, and suppose that
                $$ \{ t \in C(S) \colon t \leq s \text{ and $te$ is defined} \} \neq \emptyset. $$
            Then there exists $t \in C(S)$ such that $t = t^+s$ and $(t^+s) \bullet e = t \bullet e$ is defined. Since $S$ is associative and $(t^+s) \bullet e = (t^+s)e$, it follows that $se$ is defined. On the other hand, if $se$ is defined, then $se \leq s$ because $se = (se)^+s$ by \eqref{lr4}. Moreover, $(se) \bullet e$ is defined by Lemma \ref{lema:lr}(b) and the associativity of $S$. Therefore,
                $$ \{ t \in C(S) \colon t \leq s \text{ and $te$ is defined} \} \neq \emptyset \iff se \text{ is defined}. $$
            Suppose that $se$ is defined and let $t \in C(S)$ be such that $t \leq s$ and $t \bullet e$ is defined. By the argument above, $t^+(se)$ is defined, and in this case:
                $$ t^+(se) = (t^+s)e = te = t, $$
            where the last equality holds because $t \bullet e$ is defined, so $t = te^+ = te$. Thus, $t \leq se$, and we conclude that $s|e = se$ is the maximum element of the set $\{ t \in C(S) \colon t \leq s \text{ and $te$ is defined} \}$.\\

            \noindent\eqref{wo5} Let $e \in C(S)^+$ and suppose that $s \bullet t$ is defined. Then $st$ is defined. Since $s|e \neq \emptyset$ if and only if $se$ is defined, it follows from the associativity of $S$ that
                $$ (s \bullet t)|e \neq \emptyset \iff (st)e = (s \bullet t)e \text{ is defined} \iff te \text{ is defined} \iff t|e \neq \emptyset. $$

            \noindent\eqref{wo6} Let $e,f \in C(S)^+$ be such that $f \leq e$. Then $fe$ is defined and $f = fe$. Therefore, for each $s \in C(S)$, we have:
                $$ s|e \neq \emptyset \iff se \text{ is defined} \iff s(fe) \text{ is defined} \iff sf \text{ is defined} \iff s|f \neq \emptyset. $$

            \noindent\eqref{wo7} Let $e \in C(S)^+$ and suppose that $s \bullet t$ and $(s \bullet t)|e \neq \emptyset$. That is, $(st)e = (s \bullet t)e$ is defined. Then,
            \begin{align*}
                (x|(y|e)^+)^+ &= (x(ye)^+)^+ \\
                &= (x(ye))^+ & \eqref{lema:lr}(c) \\
                &= ((xy)e)^+ \\
                &= ((x \bullet y)|e)^+.
            \end{align*}

            \noindent\eqref{wo8} Let $e,f \in C(S)^+$ be such that $e \leq f$. Then the restriction $e|f$ is defined and $e|f = ef$. In particular, $ef$ is defined. Hence, the corestriction $e|f$ is defined and $e|f = ef$. Therefore, the restriction and corestriction coincide.\\

            \noindent\eqref{wo9} This condition is precisely \cite[Lemma 2.12(p3)]{rsgpdexpansion}. 
            
            This concludes the argument.
        \end{proof}
    \end{lemma}

    To each li-constellation we associate a left restriction semigroupoid as follows. For the li-constelation $T$, let $G(T) = T$, and define
        $$ G(T)^{(2)} = \{ (x,y) \in T \times T \colon x|y^+ \neq \emptyset \}. $$
    By \eqref{wo4}, if $x|y^+ \neq \emptyset$, then the product $(x|y^+)y^+$ is defined. Moreover, Lemma \ref{lema:constellation}(c) implies that $(x|y^+)y$ is also defined. We thus define $\otimes \colon G(T)^{(2)} \to G(T)$ by
        $$ x \otimes y = (x|y^+)y. $$
    The function $\otimes$ is referred to as the \textit{pseudo-product} of $T$. Let $+ \colon G(T) \to G(T)$ be the same function $+$ as in $T$. To differentiate the composable pairs of $T$ and $G(T)$, we write ``$xy$ is defined" when $(x,y) \in T^{(2)}$, and ``$x \otimes y$ is defined" when $(x,y) \in G(T)^{(2)}$.

    \begin{lemma} \label{prop:const-to-sgpd}
        The pair $(G(T),+)$ is a left restriction semigroupoid.

        \begin{proof}
            First we prove that $(G(T),G(T)^{(2)},\otimes)$ is a semigroupoid. Let $x,y,z \in G(T)$. Then:
            \begin{itemize}
                \item[(a)] $x \otimes y$ and $y \otimes z$ are defined if and only if $x|y^+ \neq \emptyset$ and $y|z^+ \neq \emptyset$.
                \item[(b)] $y \otimes z$ and $x \otimes (y \otimes z)$ are defined if and only if $y|z^+ \neq \emptyset$ and $x|(((y|z^+)z))^+ \neq \emptyset$.
                \item[(c)] $x \otimes y$ and $(x \otimes y) \otimes z$ are defined if and only if $x|y^+ \neq \emptyset$ and $((x|y^+)y)|z^+ \neq \emptyset$.
            \end{itemize}
            
            Suppose that $y|z^+ \neq \emptyset$, and note that $((y|z^+)z)^+ = (y|z^+)^+$ by Lemma \ref{lema:constellation}(b). From \eqref{wo4}, we have $y|z^+ \leq y$, so by \eqref{wo2}, it follows that $(y|z^+)^+ \leq y^+$. Then, applying \eqref{wo6}, we obtain
                $$ x|y^+ \neq \emptyset \iff x|((y|z^+)z)^+ = x|(y|z^+)^+ \neq \emptyset. $$
            This proves that (a) is equivalent to (b). Now suppose that $x|y^+ \neq \emptyset$. Then, using \eqref{wo5}, we get
                $$ ((x|y^+)y)|z^+ \neq \emptyset \iff y|z^+ \neq \emptyset. $$
            Hence, (a) is also equivalent to (c). Finally, let us compute:
            \begin{align*}
                x \otimes (y \otimes z) &= x \otimes (y|z^+)z \\
                &= [x|((y|z^+)z)^+][(y|z^+)z] \\
                &= [x|(y|z^+)^+][(y|z^+)z] & \ref{lema:constellation}(b) \\
                &= [(x|(y|z^+)^+)(y|z^+)]z. & \eqref{c2}
            \end{align*}
            On the other hand, recall that $(y|z^+)^+ \leq y^+$. Hence $(x|y^+)|(y|z^+)^+ = x|(y|z^+)^+$ by Proposition \ref{prop:corestriction}(c). Thus,
            \begin{align*}
                (x \otimes y) \otimes z &= (x|y^+)y \otimes z \\
                &= [((x|y^+)y)|z^+]z \\
                &= [((x|y^+)|(y|z^+)^+)(y|z^+)]z & \ref{lema:corestriction}(b) \\
                &= [(x|(y|z^+)^+)(y|z^+)]z & \ref{prop:corestriction}(c) \\
                &= x \otimes (y \otimes z).
            \end{align*}
            This proves that $G(T)$ is a semigroupoid. We now verify that $(G(T), +)$ is a left restriction semigroupoid.\\

            \noindent\eqref{lr1} For every $x \in G(T)$, we have $x^+ \in \{ t \in T \colon t \leq x^+ \text{ and $tx^+$ is defined} \}$. In particular, $x^+|x^+ \neq \emptyset$ and $x^+|x^+ = x^+$. Therefore, $x^+ \otimes x$ is defined, and
                $$ x^+ \otimes x = (x^+|x^+)x = x^+x = x, $$
            where the last equality follows from \eqref{c3}.\\

            \noindent\eqref{lr2} By \eqref{wo9}, we have that $(T^+,\leq)$ is a local meet-semillatice with the operation given by corestriciton. Therefore,
            \begin{align*}
                x^+ \otimes y^+ \text{ is defined} &\iff x^+|y^+ = x^+|(y^+)^+ \neq \emptyset \\
                &\iff y^+|(x^+)^+ = y^+|x^+ \neq \emptyset \\
                &\iff y^+ \otimes x^+ \text{ is defined}.
            \end{align*}
            In this case we have $x^+|y^+ = y^+|x^+$, $x^+|y^+ \leq x^+$, and $y^+|x^+ \leq y^+$. Moreover, since $T^+$ is closed under corestriction, it follows from Lemma \ref{lema:constellation}(a) that $x^+|y^+ = (x^+|y^+)^+$. Hence,
            \begin{align*}
                x^+ \otimes y^+ &= (x^+|y^+)y^+ \\
                &= (x^+|y^+)^+ y^+ & \ref{lema:constellation}(a) \\
                &= x^+|y^+. & \eqref{coro:restriction}
            \end{align*}
            Analogously, we have $y^+ \otimes x^+ = y^+|x^+$, and therefore $x^+ \otimes y^+ = y^+ \otimes x^+$.\\

            \noindent\eqref{lr3} Suppose that $x^+ \otimes y$ is defined. Then
            \begin{align*}
                (x^+ \otimes y)^+ &= ((x^+|y^+)y)^+ \\
                &= (x^+|y^+)^+ & \ref{lema:constellation}(b) \\
                &= x^+|y^+ & \eqref{wo9}, \ref{lema:constellation}(a) \\
                &= x^+ \otimes y^+,
            \end{align*}
            where the last equality follows from the proof of \eqref{lr2}.\\

            \noindent\eqref{lr4} Suppose that $x \otimes y$ is defined. Then $x|y^+ \neq \emptyset$ and $(x|y^+)^+ \leq x^+$. By Proposition \ref{prop:restriction}(a), we have $(x|y^+)^+ = (x|y^+)^+ | x^+$. Therefore,
            \begin{align*}
                (x \otimes y)^+ \otimes x &= ((x|y^+)y)^+ \otimes x \\
                &= (x|y^+)^+ \otimes x & \ref{lema:constellation}(b) \\
                &= ((x|y^+)^+ | x^+) x \\
                &= (x|y^+)^+ x & \ref{prop:restriction}(a) \\
                &= x|y^+ & \ref{lema:corestriction}(a) \\
                &= (x|y^+)y^+ & \eqref{c4} \\
                &= (x|(y^+)^+)y^+ & \ref{lema:constellation}(a) \\
                &= x \otimes y^+.
            \end{align*}
            This shows that $(G(T),+)$ is a left restriction semigroupoid.
        \end{proof}
    \end{lemma}

By combining Lemmas \ref{lema:sgpd-to-const} and \ref{prop:const-to-sgpd} we obtain a bijective correspondence between left restriction semigroupoids and li-constellations, described as follows.

    \begin{prop} \label{prop:objetos}
        With the above notation, we have $G(C(S)) = S$ as left restriction semigroupoids, and $C(G(T)) = (T, \leq)$ as li-constellations.

        \begin{proof}
            For a left restriction semigroupoid $S$, we have $G(C(S)) = C(S) = S$ as sets. The set of composable pairs of $G(C(S))$ is given by
            \begin{align*}
                G(C(S))^{(2)} &= \{ (s,t) \in C(S) \times C(S) \colon s|t^+ \neq \emptyset \} \\
                &= \{ (s,t) \in S \times S \colon st^+ \text{ is defined} \} \\
                &= \{ (s,t) \in S \times S \colon st \text{ is defined} \} \\
                &= S^{(2)}.
            \end{align*}
            For each pair $(s,t) \in G(C(S))^{(2)}$, the pseudo-product $s \otimes t$ is defined as
            \begin{align*}
                s \otimes t = (s|t^+) \bullet t = (st^+)t = s(t^+t) = st.
            \end{align*}
           The restriction $+ \colon G(C(S)) \to G(C(S))$ is given by the function $+ \colon C(S) \to C(S)$, which coincides with the original operation $+ \colon S \to S$. Therefore, $G(C(S)) = S$ as left restriction semigroupoids.\\

            On the other hand, for a li-constellation $(T,\leq)$, we have $C(G(T)) = G(T) = T$ as sets. The set of composable pairs of $C(G(T))$ is given by
            \begin{align*}
                C(G(T))^{(2)} &= \{ (x,y) \in G(T) \times G(T) \colon xy^+ \text{ is defined and } x \otimes y^+ = x \} \\
                &= \{ (x,y) \in T \times T \colon x|y^+ \neq \emptyset \text{ and } (x|y^+)y^+ = x \} \\
                &= \{ (x,y) \in T \times T \colon x|y^+ \neq \emptyset \text{ and } x|y^+ = x \} & \eqref{c4} \\
                &= \{ (x,y) \in T \times T \colon xy \text{ is defined} \} & \ref{prop:corestriction}(a) \\
                &= T^{(2)}.
            \end{align*}
            For each pair $(x, y) \in C(G(T))^{(2)} = T^{(2)}$, we have $x|y^+ = x$ by Proposition \ref{prop:corestriction}(a). Hence, the operation in $C(G(T))$ is given by
            \begin{align*}
                x \bullet y &= x \otimes y = (x|y^+)y = xy.
            \end{align*}
            The range function $+ \colon C(G(T)) \to C(G(T))$ is given by the restriction $+ \colon G(T) \to G(T)$, which coincides with the operation $+ \colon T \to T$. Finally, $T$ and $C(G(T))$ share the same partial order. Indeed, suppose that $x \leq y$ in $T$. Then, by Corollary \ref{coro:restriction}, $x = x^+y$, and hence
                $$ x = x^+y = x^+ \otimes y \in G(T). $$
            This shows that $x \leq y$ in $C(G(T))$. Conversely, suppose that $x \leq y$ in $C(G(T))$. Then $x^+ \otimes y$ is defined and $x = x^+ \otimes y = (x^+|y^+)y$. Since $x^+|y^+ \leq y^+$, and both $(x^+|y^+)y$ and $y^+y$ are defined, it follows from \eqref{wo1} that
                $$ x = (x^+|y^+)y \leq y^+y = y \in T. $$
            Therefore, $x \leq y$ in $T$ if and only if $x \leq y$ in $C(G(T))$. Since restrictions and corestrictions are completely determined by $\leq$, $+$ and $T^{(2)}$, it follows that $C(G(T)) = (T,\leq)$ as li-constellations.
        \end{proof}
    \end{prop}

    To establish an isomorphism between the category of left restriction semigroupoids with morphisms and the category of li-constellations with inductive radiants, it remains to establish a bijective correspondence between morphisms and inductive radiants.

    \begin{lemma} \label{lema:morphism-2}
        If $\varphi \colon S \to S'$ is a morphism, then $C\varphi = \varphi \colon C(S) \to C(S')$ is an inductive radiant.

        \begin{proof}
            \eqref{ir1} Suppose that $s \bullet t$ is defined. Then $st^+$ is defined and $st^+ = s$. From \eqref{rm1} and \eqref{rm2}, we deduce that $\varphi(s) \varphi(t)^+$ is defined and satisfies $\varphi(s)\varphi(t)^+ = \varphi(st^+) = \varphi(s)$. Hence, $\varphi(s) \bullet \varphi(t)$ is defined. Moreover, since $st^+$ is defined if and only if $st$ is defined, it follows from \eqref{rm1} that
                $$ \varphi(s) \bullet \varphi(t) = \varphi(s)\varphi(t) = \varphi(st) = \varphi(s \bullet t). $$
                
            \noindent\eqref{ir2} This conditions is precisely \eqref{rm2}.\\

            \noindent\eqref{ir3} Suppose that $s \leq t$. By Corollary \ref{coro:restriction}, we have that $s^+ \bullet t$ is defined and $s = s^+ \bullet t$. Since $s^+ \bullet t = s^+t$, it follows from \eqref{rm1} and \eqref{rm2} that
                $$ \varphi(s) = \varphi(s^+t) = \varphi(s)^+\varphi(t) \leq \varphi(t). $$

            \noindent\eqref{ir4} Suppose that $s|e \neq \emptyset$. Then $se$ is defined and, by \eqref{rm1}, it follows that $\varphi(s)\varphi(e)$ is defined and $\varphi(s)\varphi(e) = \varphi(se)$. Since $e = e^+$ by Lemma \ref{lema:lr}(b), we have $\varphi(e) = \varphi(e)^+ \in S'^+$. Hence, $\varphi(s)|\varphi(e) \neq \emptyset$, and in this case
                $$ \varphi(s|e) = \varphi(se) = \varphi(s)\varphi(e) = \varphi(s)|\varphi(e). $$
            This shows that $C\varphi = \varphi$ defines an inductive radiant from $C(S)$ to $C(S')$.
        \end{proof}
    \end{lemma}

    \begin{lemma} \label{lema:morphism-1}
        If $\varphi \colon T \to T'$ is an inductive radiant, then $G\varphi = \varphi \colon G(T) \to G(T')$ is a morphism.

        \begin{proof}
            \eqref{rm1} Suppose that $x \otimes y$ is defined. Then $x|y^+ \neq \emptyset$. By \eqref{ir4} and \eqref{ir1}, we have that $\varphi(x)|\varphi(y)^+ \neq \emptyset$ and $\varphi(x|y^+) = \varphi(x)|\varphi(y)^+$. Thus $\varphi(x) \otimes \varphi(y)$ is defined, and
            \begin{align*}
                \varphi(x \otimes y) &= \varphi((x|y^+)y) \\
                &= (\varphi(x)|\varphi(y)^+)\varphi(y) & \eqref{ir1}, \eqref{ir4}, \eqref{ir2} \\
                &= \varphi(x) \otimes \varphi(y).
            \end{align*}
            
            \noindent\eqref{rm2} This is exactly \eqref{ir2}. 
            
            Hence, $G\varphi = \varphi$ defines a morphism from $G(T)$ to $G(T')$.
        \end{proof}
    \end{lemma}

    \begin{prop} \label{prop:1st-category}
        Li-constellations (as objects) and inductive radiants (as morphisms) form a category, where the domain and range are given by the domain and range of functions, and $\circ$ is the usual composition of functions.

        \begin{proof}
            To prove that li-constellations and inductive radiants form a category, we need to show that the identity function $id_T \colon T \to T$ is an inductive radiant, and that the composition of inductive radiants is again an inductive radiant. Specifically, if $\varphi \colon T \to T'$ and $\psi \colon T' \to T''$ are inductive radiants, then their composition $\psi \circ \varphi \colon T \to T''$ is also an inductive radiant.

            It is straightforward to verify that $id_T$ is an inductive radiant. Suppose $\varphi$ and $\psi$ are inductive radiants. By Lemma \ref{lema:morphism-1}, $G\varphi$ and $G\psi$ are morphisms, and by Proposition \ref{prop:sgpd-cat-1}, the composition $G\psi \circ G\varphi$ is a morphism. Applying Lemma \ref{lema:morphism-2} and Proposition \ref{prop:objetos}, it follows that $C(G\psi \circ G\varphi) \colon T \to T''$ is an inductive radiant. By construction, 
                $$ C(G\psi \circ G\varphi)(x) = G\psi(G\varphi(x)) = \psi(\varphi(x)) = (\psi \circ \varphi)(x), \ \forall x \in T. $$
            Therefore, $\psi \circ \varphi = C(G\psi \circ G\varphi)$ is an inductive radiant.
        \end{proof}
    \end{prop}

   To conclude this section we present below the first part of our ESN-type theorem.

    \begin{theorem} \label{teo:1st-isomorphism}
        The category of left restriction semigroupoids and morphisms is isomorphic to the category of li-constellations and inductive radiants.

        \begin{proof}
            From Proposition \ref{prop:objetos}, Lemma \ref{lema:morphism-2}, and Lemma \ref{lema:morphism-1}, we know that $C$ and $G$ are mutually inverse bijections between left restriction semigroupoids and li-constellations, and between morphisms and inductive radiants. Therefore, to establish the categorical isomorphism, it remains to show that $C$ and $G$ are functors.

            We have already observed that $Gid_T = id_{G(T)}$, and that $C(G\psi \circ G\varphi) = \psi \circ \varphi$ whenever the composition $\psi \circ \varphi$ is defined. Applying $G$ to this identity and using the fact that $GC\rho = \rho$ for every inductive radiant $\rho$, we obtain:
                $$ G(\psi \circ \varphi) = GC(G\psi \circ G\varphi) = G\psi \circ G\varphi. $$
         This shows that $G$ is a functor. Similarly, one can show that $C$ is also a functor, completing the proof of the categorical isomorphism.
        \end{proof}
    \end{theorem}
    
\section{The Szendrei Expansion of li-constellations}
In this section we recall the construction of the Szendrei expansion of left restriction semigroupoids, as presented in \cite{rsgpdexpansion}. We then introduce the Szendrei expansion of li-constellations and establish a universal property that characterizes this expansion. This result is included in the present work because it plays a key role in the proof of the second part of our ESN-type theorem.\\

    The Szendrei expansion of a left restriction semigroupoid $S$ is constructed as follows. For each $e \in S^+$, let $\mathcal{P}_e(S)$ denote the family of finite subsets $A$ of $S$ such that $e \in A$ and $a^+ = e$ for all $a \in A$. Define:
        $$ Sz(S) = \left\{ (A,a) \colon A \in \bigcup_{e \in S^+} \mathcal{P}_e(S), a \in A \right\}, $$
    equipped with the following partial binary operation:
        $$ (A,a)(B,b) = \begin{cases}
            ((ab)^+ A \cup aB, ab), & \text{if $ab$ is defined}, \\
            \text{undefined}, & \text{otherwise}.
        \end{cases} $$
    With this operation, $Sz(S)$ forms a semigroupoid. Moreover, if we define $(A, a)^+ = (A, a^+)$ for each $(A, a) \in Sz(S)$, then $(Sz(S), +)$ becomes a left restriction semigroupoid \cite[Proposition 4.6]{rsgpdexpansion}.\\

    Since the Szendrei expansion $Sz(S)$ of the left restriction semigroupoid $S$ is itself a left restriction semigroupoid, we may apply Proposition \ref{prop:objetos} to obtain the corresponding li-constellation $C(Sz(S))$. We consider here the same functors $C$ and $G$ introduced in Section 4. In this setting, the structure of $C(Sz(S))$ can be described in terms of the structure of $C(S)$. Indeed, we have:
    \begin{align*}
        &\quad (A,a) \bullet (B,b) \text{ is defined} \\
        &\iff (A,a)(B,b) \text{ is defined and } (A,a)(B,b)^+ = (A,a) \\
        &\iff ab \text{ is defined and } ((ab^+)^+ A \cup ab, ab^+) = (A,a) \\
        &\iff ab \text{ is defined, } ab^+ = a \text{ and } aB \subseteq A \\
        &\iff a \bullet b \text{ is defined and } a \bullet B \subseteq A.
    \end{align*}
    In this case,
        $$ (A,a) \bullet (B,b) = ((ab)^+ A \cup aB, ab) = (A,a \bullet b). $$

        The range function $+ \colon C(Sz(S)) \to C(Sz(S))$ is the restriction of $+ \colon Sz(S) \to Sz(S)$, induced by the restriction $+ \colon S \to S$, which also defines the range function on $C(S)$. The partial order in $C(Sz(S))$ is defined by:
    \begin{align*}
        (A,a) \leq (B,b) &\iff (A,a)^+(B,b) \text{ is defined and } (A,a) = (A,a)^+(B,b) \\
        &\iff a^+b \text{ is defined and } (A,a) = ((a^+b)^+ A \cup a^+B, a^+b) \\
        &\iff a^+b \text{ is defined, } a^+b = a \text{ and } a^+B \subseteq A \\
        &\iff a \leq b \text{ and } a^+ \bullet B \subseteq A.
    \end{align*}
    Let $(E,e) \in Sz(S)^+$, and suppose that $(E,e) \leq (A,a)^+$. Then $e \in S^+$ and $e \leq a^+$. In particular, $ea$ is defined. Thus, the restriction of $(A,a)$ to $(E,e)$ is given by:
        $$ (E,e)|(A,a) = (E,e)(A,a) = ((ea)^+ E \cup eA, ea) = (E,e|a). $$
    Moreover, we have $(A,a)|(E,e) \neq \emptyset$ if and only if $ae$ is defined. In that case, we have
        $$ (A,a)|(E,e) = (A,a)(E,e) = ((ae)^+ A \cup aE,ae) = ((a|e)^+ A \cup (a|e)E, a|e). $$
    This proves the following result.

    \begin{prop} \label{prop:li-constellation-expansion}
        Given a li-constellation and an element $e \in T^+$, denote by $\mathcal{P}_e(T)$ the collection of finite subsets $A \subseteq T$ such that $e \in A$ and $a^+ = e$ for every $a \in A$. Define
            $$ Sz(T) = \left\{ (A,a) \colon A \in \bigcup_{e \in T^+} \mathcal{P}_e(T), a \in A \right\}, $$
        equipped with the following partial binary operation
            $$ (A,a)(B,b) = \begin{cases}
                (A,ab), & \text{if $ab$ is defined}, \\
                \text{undefined}, & \text{otherwise}.
            \end{cases} $$
        Define $(A,a)^+ = (A,a^+)$, for every $(A,a) \in Sz(T)$. Then $Sz(T)$ is a left constellation. Define a partial order $\leq$ in $Sz(T)$ by
            $$ (A,a) \leq (B,b) \iff a \leq b \text{ and } a^+B \subseteq A. $$
        With this order, $(Sz(T), \leq)$ is a li-constellation.
    \end{prop}

    The li-constellation $(Sz(T),\leq)$ constructed in Proposition~\ref{prop:li-constellation-expansion} is called the \textit{Szendrei expansion} of $(T,\leq)$. By construction,  it satisfies $Sz(T) = Sz(CG(T)) = C(Sz(G(T)))$. Therefore, restrictions and corestrictions in $Sz(T)$ are given by
        $$ (E,e)|(A,a) = (E,e|a) \quad\text{and}\quad (A,a)|(E,e) = ((a|e)^+A \cup (a|e)E, a|e), $$
    when defined.

    \begin{lemma} \label{lema:iota}
        Let $\iota \colon T \to Sz(T)$ be defined by $\iota(x) = (\{x^+,x\},x)$ for every $x \in T$. Then $\iota$ is an inductive preradiant. Moreover, whenever $x|y^+ \neq \emptyset$, we have
            $$ \iota(x) \otimes \iota(y) = \iota(x)^+ \otimes \iota(x \otimes y), $$
        where $\otimes$ denotes the composition in $G(T)$.

        \begin{proof}
            \eqref{ip1} Suppose that $xy$ is defined and recall that $(xy)^+ = x^+$ by Lemma \ref{lema:constellation}(b). Then
            \begin{align*}
                \iota(x) &= (\{x^+,x\},x), & \iota(x)^+ &= (\{x^+,x\},x^+), \\
                \iota(y) &= (\{y^+,y\},y), & \iota(y)^+ &= (\{y^+,y\},y^+), \\
                \iota(xy) &= (\{x^+,xy\},xy), & \iota(xy)^+ &= (\{x^+,xy\},x^+).
            \end{align*}
            Since $xy$ is defined, it follows from Proposition \ref{prop:corestriction}(a) that $x|y^+ \neq \emptyset$ and $x|y^+ = x$. Therefore, $\iota(x)|\iota(y)^+ \neq \emptyset$, and we compute:
            \begin{align*}
                \iota(x)|\iota(y)^+ &= (\{x^+,x\},x)|(\{y^+,y\},y^+) \\
                &= ( \{(x|y^+)^+x^+, (x|y^+)^+x, (x|y^+)y^+ ,(x|y^+)y\} , x|y^+) \\
                &= ( \{x^+x^+, x^+x, xy^+ ,xy\} , x) \\
                &= ( \{x^+,x,xy\}, x).
            \end{align*}
            In the last step, we use that $x^+x = x$ by \eqref{c3}, $x^+x^+ = x^+$ and $xy^+ = x$ by \eqref{c4}. On the other hand, since $x^+|x^+ \neq \emptyset$ and $x^+|x^+ = x^+$, we have $\iota(x)^+|\iota(xy)^+ \neq \emptyset$, and:
            \begin{align*}
                \iota(x)^+|\iota(xy)^+ &= (\{x^+,x\},x^+)|(\{x^+,xy\},x^+) \\
                &= ( (x^+|x^+)^+\{x^+,x\} \cup (x^+|x^+)\{x^+,xy\}, x^+|x^+ ) \\
                &= ( \{ x^+x^+, x^+x , x^+(xy) \}, x^+ ) \\
                &= ( \{ x^+, x, xy \}, x^+ ),
            \end{align*}
            where in the last equality we use \eqref{c2} and \eqref{c3} to conclude that $x^+(xy) = (x^+x)y = xy$. Hence,
            \begin{align*}
                \iota(x)^+ \otimes \iota(xy) &= (\iota(x)^+|\iota(xy)^+)\iota(xy) \\
                &= (\{x^+, x, xy\}, x^+)(\{x^+, xy\}, xy) \\
                &= ( \{x^+, x, xy\}, x^+(xy) ) \\
                &= ( \{x^+, x, xy\}, xy ) \\
                &= ( \{x^+, x, xy\}, x )( \{y^+,y\}, y ) \\
                &= (\iota(x)|\iota(y)^+)\iota(y) \\
                &= \iota(x) \otimes \iota(y).
            \end{align*}

            \noindent\eqref{ip2} Note that $\iota(x)^+ = ({x^+, x}, x^+)$ and $\iota(x^+) = ({x^+}, x^+)$. Since $x^+ \leq x^+$ and $x^+ {x^+} = {x^+} \subseteq {x^+, x}$, it follows from the definition of the order on $Sz(T)$ that $\iota(x)^+ \leq \iota(x^+)$.\\

            \noindent\eqref{ip3} Suppose that $x \leq y$. Then, by \eqref{wo2}, we have $x^+ \leq y^+$, and by Corollary \ref{coro:restriction}, it follows that $x = x^+y$ and $x^+ = x^+y^+$. Hence, $x^+\{y^+,y\} = \{x^+y^+,x^+y\} = \{x^+,x\}$, which shows that
                $$ \iota(x) = (\{x^+,x\},x) \leq (\{y^+,y\},y) = \iota(y). $$

            \noindent\eqref{ip4} Suppose that $x|e \neq \emptyset$. Then $\iota(x)|\iota(e) \neq \emptyset$, and we compute: 
            \begin{align*}
                \iota(x)|\iota(e) &= (\{x^+,x\},x)|(\{e\},e) \\
                &= ( (x|e)^+\{x^+,x\} \cup (x|e)\{e\}, x|e ) \\
                &= (\{ (x|e)^+x^+, (x|e)^+x, (x|e)e \}, x|e ) \\
                &= (\{ (x|e)^+, (x|e)^+x, x|e \}, x|e ) & \eqref{c4} \\
                &= (\{ (x|e)^+, x|e \}, x|e ) \\
                &= \iota(x|e).
            \end{align*}
            In the fifth equality we use the fact that $x|e \leq x$, and therefore $(x|e)^+x = x|e$ by Corollary \ref{coro:restriction}.\\

            \noindent\eqref{ip5} Let $e \leq x^+$. Then $e = e^+|x^+$ by Proposition \ref{prop:restriction}(a). Since $e^+,x^+ \in T^+$ and $e^+|x^+ \neq \emptyset$, it follows from \eqref{wo9} that $e^+|x^+ \in T^+$. Applying Lemma \ref{lema:constellation}(a) to both $e^+|x^+$ and to $e$, we obtain
                $$ e = e^+|x^+ = (e^+|x^+)^+ = (e|x^+). $$
            On the other hand, we have $e|x = ex$ by Proposition \ref{prop:restriction}(d). Therefore,
            \begin{align*}
                \iota(e) | \iota(x)^+ &= (\{e\},e)|(\{x^+,x\},x^+) \\
                &= ( \{ (e|x^+)^+e , (e|x^+)x^+, (e|x^+)x \}, e|x^+ ) \\
                &= ( \{ ee, ex^+, ex \}, e ) \\
                &= ( \{ e, e|x \}, e ) \\
                &= ( \{(e|x)^+, e|x\}, (e|x)^+ ) \\
                &= \iota(e|x)^+.
            \end{align*}
            This shows that $\iota \colon T \to Sz(T)$ is an inductive preradiant.\\
            
            Suppose that $x,y \in T$ are such that $x|y^+ \neq \emptyset$. Since $x|y^+ \leq x$, it follows from \eqref{wo2} that $(x|y^+)^+ \leq x^+$. By Corollary \ref{coro:restriction}, we have
                $$ (x|y^+)^+ = (x|y^+)^+x^+. $$
            Moreover, from Proposition \ref{prop:restriction}(a) and \eqref{wo9}, we obtain
                $$ x^+|(x|y^+)^+ = (x|y^+)^+|x^+ = (x|y^+)^+. $$
            In particular, $x^+|(x|y^+)^+ \neq \emptyset$. Since $(x \otimes y)^+ = ((x|y^+)y)^+ = (x|y^+)^+$ by Lemma \ref{lema:constellation}(a), it follows that $\iota(x)^+ | \iota(x \otimes y)^+ \neq \emptyset$. In this case,
            \begin{align*}
                &\quad \iota(x)^+ | \iota(x \otimes y)^+ \\
                &= (\{x^+,x\},x^+)|(\{(x|y^+)^+, (x|y^+)y\}, (x|y^+)^+) \\
                &= ( (x^+|(x|y^+)^+)^+\{  x^+, x\} \cup (x^+|(x|y^+)^+) \{(x|y^+)^+, (x|y^+)y \}, x^+|(x|y^+)^+ ) \\
                &= ( (x|y^+)^+\{  x^+, x\} \cup (x|y^+)^+ \{(x|y^+)^+, (x|y^+)y \}, (x|y^+)^+ ) \\
                &= ( \{ (x|y^+)^+, (x|y^+)^+x, (x|y^+)y \}, (x|y^+)^+ ).
            \end{align*}
            Therefore, we have
            \begin{align*}
                &\quad \iota(x)^+ \otimes \iota(x \otimes y) \\
                &= (\iota(x)^+|\iota(x \otimes y)^+) \iota(x \otimes y) \\
                &= ( \{ (x|y^+)^+, (x|y^+)^+x, (x|y^+)y \}, (x|y^+)^+ )( \{ (x|y^+)^+, (x|y^+)y \}, (x|y^+)y ) \\
                &= ( \{ (x|y^+)^+, (x|y^+)^+x, (x|y^+)y \}, (x|y^+)^+ (x|y^+)y ) \\
                &= ( \{ (x|y^+)^+, (x|y^+)^+x, (x|y^+)y \}, (x|y^+)y ).
            \end{align*}
            On the other hand,
            \begin{align*}
                &= \iota(x) \otimes \iota(y) \\
                &= (\iota(x)|\iota(y)^+)\iota(y) \\
                &= ( \{(x|y^+)^+ x^+, (x|y^+)^+ x, (x|y^+) y^+, (x|y^+)y\}, x|y^+ )(\{y^+,y\},y) \\
                &= ( \{(x|y^+)^+, (x|y^+)^+ x, (x|y^+)y\}, x|y^+ )(\{y^+,y\},y) \\
                &= ( \{(x|y^+)^+, (x|y^+)^+ x, (x|y^+)y\}, (x|y^+)y ).
            \end{align*}
            We thus conclude that $\iota(x) \otimes \iota(y) = \iota(x)^+ \otimes \iota(x \otimes y)$ whenever $x|y^+ \neq \emptyset$.
        \end{proof}
    \end{lemma}

    \begin{lemma} \label{lema:generated}
        The Szendrei expansion $Sz(T)$ is generated by $\iota(T)$ under composition, $+$ and corestriction.
        
        \begin{proof}
            Let $(A,a) \in Sz(T)$. Then $\iota(a) = (\{a^+,a\},a)$ and $a^+\{a^+,a\} = \{a^+,a\} \subseteq A$. Hence, the composition $(A,a^+)\iota(a)$ is defined and
                $$ (A,a)^+ \iota(a) = (A,a^+)(\{a^+,a\},a) = (A,a^+a) = (A,a). $$
            On the other hand, let $(E,e), (F,e) \in Sz(T)$ for some $e \in T^+$. Since $e|e \neq \emptyset$, we have $(E,e)|(F,e) \neq \emptyset$. In this case,
                $$ (E,e)|(F,e) = ( (e|e)^+ E \cup (e|e)F, e|e ) = (E \cup F,e). $$
            Let $A = \{ a^+, a_1, \dots, a_n\}$ with $a_j^+ = a^+$. Then
            \begin{gather*}
                ((\iota(a_1)^+ | \iota(a_2)^+) | \dots )| \iota(a_n)^+ = \left( \bigcup_{j=1}^n \{ a^+, a_j \}, a^+ \right) = (A,a^+) = (A,a)^+.
            \end{gather*}
            This shows that $(A,a) = ( \iota(a_1)^+|\dots|\iota(a_n)^+ )\iota(a)$ can be constructed from $\iota(T)$ using $+$, corestriction, and composition.
        \end{proof}
    \end{lemma}

    \begin{lemma} \label{lema:composition}
        If $\varphi \colon T \to T'$ is an inductive preradiant and $\Phi \colon T' \to T''$ is an inductive radiant, then $\Phi \varphi \colon T \to T''$ is an inductive preradiant.

        \begin{proof}
            \eqref{ip1} Suppose that $xy$ is defined. Then $\varphi(x) | \varphi(y)^+ \neq \emptyset$ by \eqref{ip1}, and $$\Phi\varphi(x) | \Phi\varphi(y)^+ \neq \emptyset$$ by \eqref{ir4}. Similarly, we have $\Phi\varphi(x)^+ | \Phi\varphi(xy)^+ \neq \emptyset$. Therefore, both $\Phi\varphi(x)^+ \otimes \Phi\varphi(xy)$ and $\Phi\varphi(x) \otimes \Phi\varphi(y)$ are defined. Since $\Phi$ preserves composition \eqref{ir1}, the range operation \eqref{ir2}, and corestriction \eqref{ir4}, it follows that
            \begin{align*}
                \Phi\varphi(x)^+ \otimes \Phi\varphi(xy) &= ( \Phi\varphi(x)^+ | \Phi\varphi(xy)^+ ) \Phi\varphi(xy) \\
                &= \Phi( (\varphi(x)^+|\varphi(xy)^+)\varphi(xy) ) \\
                &= \Phi(\varphi(x)^+ \otimes \varphi(xy)) \\
                &= \Phi(\varphi(x) \otimes \varphi(y)) & \eqref{ip1} \\
                &= \Phi((\varphi(x)|\varphi(y)^+)\varphi(y)) \\
                &= (\Phi\varphi(x)|\Phi\varphi(y)^+)\Phi\varphi(y) \\
                &= \Phi\varphi(x) \otimes \Phi\varphi(y).
            \end{align*}

            \noindent\eqref{ip2} From \eqref{ip2} applied to $\varphi$, we have $\varphi(x)^+ \leq \varphi(x^+)$. Then, by \eqref{ir2} and \eqref{ir3}, it follows that $\Phi\varphi(x)^+ \leq \Phi\varphi(x^+)$.\\

            \noindent\eqref{ip3} If $x \leq y$, then $\varphi(x) \leq \varphi(y)$ by \eqref{ip3} for $\varphi$. Therefore, by \eqref{ir3}, it follows that $\Phi\varphi(x) \leq \Phi\varphi(y)$.\\

            \noindent\eqref{ip4} Suppose that $x|e \neq \emptyset$. Then $\varphi(x)|\varphi(e) \neq \emptyset$ by \eqref{ip4}, and $\Phi\varphi(x)|\Phi\varphi(e) \neq \emptyset$ by \eqref{ir4}. Moreover, we have $\Phi\varphi(x|e) = \Phi(\varphi(x)|\varphi(e)) = \Phi\varphi(x)|\Phi\varphi(e)$.\\

            \noindent\eqref{ip5} Let $e \leq x^+$. Then the restriction $e|x$ is defined, and
            \begin{align*}
                \Phi\varphi(e|x)^+ &= \Phi(\varphi(e|x)^+) & \eqref{ir2} \\
                &= \Phi(\varphi(e)|\varphi(x)^+) & \eqref{ip5} \\
                &= \Phi\varphi(e)|\Phi(\varphi(x)^+). & \eqref{ir4}, \eqref{ir2}
            \end{align*}
            This shows that $\Phi\varphi$ is an inductive preradiant.
        \end{proof}
    \end{lemma}

    In what follows we present the main result of this section. It establishes a universal property that characterizes the Szendrei expansion $Sz(T)$. This result generalizes both \cite[Theorem 5.1]{gould2011actions} and \cite[Proposition 4.6]{gilbert2005actions}, which address similar universal properties in the context of inductive constellations and inductive groupoids, respectively.

    \begin{theorem} \label{teo:orderedSzendrei}
        For every inductive preradiant $\varphi \colon T \to T'$, there exists a unique inductive radiant $\Phi \colon Sz(T) \to T'$ such that $\Phi \iota = \varphi$. That is, the following diagram commutes:
        \begin{center}
            \begin{tikzpicture}
                \tikzstyle{every path}=[draw,->];

                \node (T) at (0,0) {$T$};
                \node (T2) at (2,0) {$T'$};
                \node (ST) at (0,-2) {$Sz(T)$};

                \path (T) to node[left]{$\iota$} (ST);
                \path (T) to node[above]{$\varphi$} (T2);
                \path[dashed] (ST) to node[below right]{$\Phi$} (T2);
            \end{tikzpicture}
        \end{center}
        Conversely, for every inductive radiant $\Phi \colon Sz(T) \to T'$, the composition $\Phi \iota \colon T \to T'$ is an inductive preradiant.

        \begin{proof}
            Let $\varphi \colon T \to T'$ be an inductive preradiant. We define a function $\Phi \colon Sz(T) \to T'$ by
                $$ \Phi(A,x) = \left( \wedge_{a \in A} \varphi(a)^+ \right) \varphi(x), $$
            where $\wedge$ denotes the corestriction in $(T',\leq)$. We begin by showing that $\Phi$ is well defined.
            
            Indeed, for every $a \in A$, we have $a^+ = x^+$, hence $f(a)^+ \leq f(a^+) = f(x^+)$ by \eqref{ip2}. That is, $\varphi(a)^+$ belongs to the same connected component of $T'^+$ as $\varphi(x^+)$, and thus the corestriction $\varphi(a)^+ \wedge \varphi(b)^+$ is defined for all $a,b \in A$. Since $(T'^+,\leq)$ is a local meet-semilattice, we have $\varphi(a)^+ \wedge \varphi(b)^+ \in T'^+$. Moreover, as $A$ is finite,  it follows by induction that $\wedge_{a \in A} \varphi(a)^+$ is defined. Since $\wedge_{a \in A} \varphi(a)^+ \leq \varphi(a)^+$ for all $a \in A$ and $x \in A$, we conclude that $\wedge_{a \in A} \varphi(a)^+ \leq \varphi(x)^+$. Then, by \eqref{wo3}, the restriction $\left( \wedge_{a \in A} \varphi(a)^+ \right)|\varphi(x)$ is defined, and by Proposition \ref{prop:restriction}(d) we obtain
                $$ \left( \wedge_{a \in A} \varphi(a)^+ \right)|\varphi(x) = \left( \wedge_{a \in A} \varphi(a)^+ \right)\varphi(x) = \Phi(A,x). $$
            Hence, $\Phi$ is a well defined function. Hence, $\Phi$ is a well-defined function. 
            
            We now proceed to show that $\Phi$ is an inductive radiant. The proof of \eqref{ir1} will be deferred to the end of the argument.\\

            \noindent\eqref{ir2} Let $(A,e) \in Sz(T)$ with $e \in T^+$. Since $e = e^+$ by Lemma \ref{lema:constellation}(a), it follows from \eqref{ip2} that $\varphi(e) \geq \varphi(e)^+$.  From the previous argument we also have $\varphi(e) \geq \varphi(e)^+ \geq \wedge_{a \in A} \varphi(a)^+$, and thus
                $$ \Phi(A,e) = \left( \wedge_{a \in A} \varphi(a)^+ \right) \varphi(e) = \left( \wedge_{a \in A} \varphi(a)^+ \right) \wedge \varphi(e) = \wedge_{a \in A} \varphi(a)^+ \in T'^+. $$
            Now, for $(A,x) \in Sz(T)$, we have $\Phi(A,x) = \left( \wedge_{a \in A} \varphi(a)^+ \right) \varphi(x) = \Phi(A,x^+) \varphi(x)$. Hence,
                $$ \Phi(A,x)^+ = \left( \Phi(A,x^+)\varphi(x) \right)^+ = \Phi(A,x^+)^+ = \Phi(A,x^+), $$
            where the second equality follows from Lemma \ref{lema:constellation}(b), and the last from Lemma \ref{lema:constellation}(a), together with the fact that $\Phi(A, x^+) \in T'^+$.\\

            \noindent\eqref{ir3} Let $(A, x) \leq (B, y)$. Then $x \leq y$ and $a^+ B \subseteq A$. By \eqref{wo2}, we have $x^+ \leq y^+ = b^+$ for every $b \in B$. Therefore, the restriction $x^+ | b$ is defined and equals $x^+ b$ by Proposition \ref{prop:restriction}(d). Since $x^+ b = x^+ | b \leq b$, it follows from \eqref{ip3} and \eqref{wo2} that
                $$ \varphi(x^+b) \leq \varphi(b) \quad\text{and}\quad \varphi(x^+b)^+ \leq \varphi(b)^+. $$
            On the other hand, we have $\varphi(x) \leq \varphi(y)$ by \eqref{ip3}. Hence, by Proposition \ref{prop:restriction}(a), it follows that $\varphi(x) = \varphi(x)^+ | \varphi(y)$. Now, from Proposition \ref{prop:restriction}(c), we obtain that whenever $e \in T'^+$ and $e \leq \varphi(x)^+ \leq \varphi(y)^+$, we have
                $$ e|\varphi(x) = e|(\varphi(x)^+|\varphi(y)) = e|\varphi(y). $$
            Since $\wedge_{a \in A} \varphi(a)^+ \in T'^+$ and $\wedge_{a \in A} \varphi(a)^+ \leq \varphi(x)^+ \leq \varphi(y)^+$, we conclude that
            \begin{align*}
                \Phi(A,x) &= \left( \wedge_{a \in A} \varphi(a)^+ \right) | \varphi(x) \\
                &= \left( \wedge_{a \in A} \varphi(a)^+ \right) | \varphi(y) \\
                &\leq \left( \wedge_{b \in B} \varphi(x^+b)^+ \right) | \varphi(y) & (x^+B \subseteq A) \\
                &\leq \left( \wedge_{b \in B} \varphi(b)^+ \right) | \varphi(y) \\
                &= \Phi(B,y).
            \end{align*}

            \noindent\eqref{ir4} If $(A,x)|(B,e) \neq \emptyset$, then $x|e \neq \emptyset$ and $(A,x)|(B,e) = ( (x|e)^+A \cup (x|e)B, x|e )$. Therefore,
            \begin{align}
                \Phi((A,x)|(B,e)) = \left( \wedge_{a \in A} \varphi( (x|e)^+a )^+ \right) \wedge \left( \wedge_{b \in B} \varphi((x|e)b)^+ \right) \varphi(x|e). \label{eq:Fcr}
            \end{align}
           On the other hand, from $x|e \neq \emptyset$ we have $\varphi(x) | \varphi(e) \neq \emptyset$ by \eqref{ip4}. Since $\Phi(B, e) \in T'^+$ by the previous argument, $\varphi(e) \in T'^+$ by Proposition \ref{prop:preradiant}, and $\Phi(B, e) = \wedge_{b \in B} \varphi(b)^+ \leq \varphi(e)$, it follows from \eqref{wo6} that $\varphi(x) | \Phi(B, e) \neq \emptyset$. Moreover, since $\Phi(A, x) = \Phi(A, x^+), \varphi(x)$, it follows from \eqref{wo5} that $\Phi(A, x) | \Phi(B, e) \neq \emptyset$, and in this case
            \begin{align*}
                \Phi(A,x)|\Phi(B,e) &= \left( \left( \wedge_{a \in A} \varphi(a)^+ \right) \varphi(x) \right) | \left( \left( \wedge_{b \in B} \varphi(b)^+ \right) \varphi(e) \right).
            \end{align*}
           Now we use the fact that composition and corestriction in $T'$ are defined via the pseudo-product $\otimes$, that is, as the composition in the left restriction semigroupoid $G(T')$. Therefore, we can rewrite $\Phi(A, x) | \Phi(B, e)$ as
            \begin{align*}
                \Phi(A,x)|\Phi(B,e) &= \left( \left( \wedge_{a \in A} \varphi(a)^+ \right) \varphi(x) \right) | \left( \left( \wedge_{b \in B} \varphi(b)^+ \right) \varphi(e) \right) \\
                &= \left( \otimes_{a \in A} \varphi(a)^+ \right) \otimes \varphi(x) \otimes \left( \otimes_{b \in B} \varphi(b)^+ \right) \otimes \varphi(e) \\
                &= \left( \otimes_{a \in A} \varphi(a)^+ \right) \otimes \left( \otimes_{b \in B} (\varphi(x) \otimes \varphi(b))^+ \right) \otimes \varphi(x) \otimes \varphi(e). & \eqref{lr4}
            \end{align*}
    We observe that in the last equality, it is necessary to apply \eqref{lr4} as many times as there are elements in $B$. On the other hand, since $\varphi((x|e)e)^+ = \varphi(x|e)^+$ by \eqref{c4}, and the meet operation satisfies $a \wedge (a \wedge b) = a \wedge b$, it follows that
            \begin{align*}
                \wedge_{b \in B} \varphi((x|e)b)^+ &= \wedge_{b \in B} (\varphi(x|e)^+ \wedge \varphi((x|e)b)^+) \\
                &= \otimes_{b \in B} (\varphi(x|e)^+ \otimes \varphi((x|e)b)^+) \\
                &= \otimes_{b \in B} (\varphi(x|e)^+ \otimes \varphi((x|e)b))^+ & \eqref{lr3} \\
                &= \otimes_{b \in B} (\varphi(x|e) \otimes \varphi(b))^+ & \eqref{ip1} \\
                &= \otimes_{b \in B} (\varphi(x)|\varphi(e) \otimes \varphi(b))^+ & \eqref{ip4} \\
                &= \otimes_{b \in B} (\varphi(x) \otimes \varphi(e) \otimes \varphi(b)^+)^+ & \ref{lema:lr}(c) \\
                &= \otimes_{b \in B} (\varphi(x) \otimes \varphi(b)^+)^+ \\
                &= \otimes_{b \in B} (\varphi(x) \otimes \varphi(b))^+, & \ref{lema:lr}(c)
            \end{align*}
            In the second-to-last equality we use \eqref{ip2} to obtain $\varphi(b)^+ \leq \varphi(b^+) = \varphi(e)$, and hence $\varphi(e) \otimes \varphi(b)^+ = \varphi(b)^+$ in the left restriction semigroupoid $(G(T'), \otimes)$. 
            
            Since $(x|e)^+ \leq x^+$ by \eqref{wo2}, and $x^+ = a^+$ for every $a \in A$, the restriction $(x|e)^+ | a$ is defined, and $(x|e)^+ | a = (x|e)^+ a$ by Proposition \ref{prop:restriction}(d). Therefore, it follows from \eqref{ip5} that
                $$ \varphi((x|e)^+a)^+ = \varphi((x|e)^+ | a)^+ = \varphi((x|e)^+) | \varphi(a)^+ = \varphi((x|e)^+) \otimes \varphi(a)^+. $$
          By substituting these identities into \eqref{eq:Fcr}, we obtain
            \begin{align*}
                &\quad \Phi((A,x)|(B,e)) \\
                &= \left( \wedge_{a \in A} \varphi( (x|e)^+a )^+ \right) \wedge \left( \wedge_{b \in B} \varphi((x|e)b)^+ \right) \varphi(x|e) \\
                &= \left( \otimes_{a \in A} (\varphi((x|e)^+) \otimes \varphi(a)^+) \right) \otimes \left( \otimes_{b \in B} (\varphi(x) \otimes \varphi(b))^+ \right) \otimes \varphi(x) \otimes \varphi(e) \\
                &= \left( \otimes_{a \in A} \varphi(a)^+ \right) \otimes \varphi((x|e)^+) \otimes \left( \otimes_{b \in B} (\varphi(x) \otimes \varphi(b))^+ \right) \otimes \varphi(x) \otimes \varphi(e). & \eqref{lr2}
            \end{align*}
            From \eqref{ip2}, the definition of meet and the previous calculations, we have
                $$ \varphi((x|e)^+) \geq \varphi(x|e)^+ \geq \varphi(x|e)^+ \wedge \varphi((x|e)b)^+ = (\varphi(x) \otimes \varphi(b))^+. $$
            Hence, we can omit the term $\varphi((x|e)^+)$ in the expression for $\Phi((A, x) | (B, e))$, and conclude that $\Phi((A,x)|(B,e)) = \Phi(A,x)|\Phi(B,e)$.\\

            \noindent\eqref{ip1} Suppose that $(A, x)(B, y)$ is defined. Then $xy$ is defined and $xB \subseteq A$. By Proposition \ref{prop:corestriction}(a), we have that $\Phi(A, x), \Phi(B, y)$ is defined if and only if
                $$ \Phi(A,x)|\Phi(B,y)^+ \neq \emptyset \text{ and } \Phi(A,x)|\Phi(B,y)^+ = \Phi(A,x). $$
           Similarly, we have that $(A, x) | (B, y)^+ \neq \emptyset$ and $(A, x) = (A, x) | (B, y)^+$. Since $\Phi$ preserves corestriction \eqref{ir4} and range \eqref{ir2}, it follows that
                $$ \Phi(A,x) = \Phi((A,x)|(B,y)^+) = \Phi(A,x)|\Phi((B,y)^+) = \Phi(A,x)|\Phi(B,y)^+. $$
            Therefore, $\Phi(A, x), \Phi(B, y)$ is defined. It remains to show that $\Phi((A,x)(B,y)) = \Phi(A,x)\Phi(B,y)$. Using the fact that both composition and corestriction in $T'$ are given by the pseudo-product $\otimes$, we obtain:
            \begin{align*}
                &\quad \Phi(A,x)\Phi(B,y) \\
                &= \left( \left( \wedge_{a \in A} \varphi(a)^+ \right) \varphi(x) \right)\left( \left( \wedge_{b \in B} \varphi(b)^+ \right) \varphi(y) \right) \\
                &= \left( \otimes_{a \in A} \varphi(a)^+ \right) \otimes \varphi(x) \otimes \left( \otimes_{b \in B} \varphi(b)^+ \right) \otimes \varphi(y) \\
                &= \left( \otimes_{a \in A} \varphi(a)^+ \right) \otimes \left( \otimes_{b \in B} (\varphi(x) \otimes \varphi(b))^+ \right) \otimes \varphi(x) \otimes \varphi(y) & \eqref{lr4} \\
                &= \left( \otimes_{a \in A} \varphi(a)^+ \right) \otimes \left( \otimes_{b \in B} (\varphi(x)^+ \otimes \varphi(xb))^+ \right) \otimes \varphi(x)^+ \otimes \varphi(xy) & \eqref{ip1} \\
                &= \left( \otimes_{a \in A} \varphi(a)^+ \right) \otimes \left( \otimes_{b \in B} (\varphi(x)^+ \otimes \varphi(xb)^+) \right) \otimes \varphi(x)^+ \otimes \varphi(xy) & \eqref{lr3} \\
                &= \left( \otimes_{a \in A} \varphi(a)^+ \right) \otimes \varphi(x)^+ \otimes \left( \otimes_{b \in B} \varphi(xb)^+ \right) \otimes \varphi(xy) & \eqref{lr2}, \ref{lema:lr}(a) \\
                &= \left( \otimes_{a \in A} \varphi(a)^+ \right) \otimes \varphi(xy).
            \end{align*}
           In the last identity we use the fact that $x \in A$ and $xB \subseteq A$. Since composition in $Sz(T)$ is defined by $(A, x)(B, y) = (A, xy)$, we conclude that
                $$ \Phi((A,x)(B,y)) = \left( \wedge_{a \in A} \varphi(a)^+ \right) \varphi(xy) = \left( \otimes_{a \in A} \varphi(a)^+ \right) \otimes \varphi(xy) = \Phi(A,x)\Phi(B,y). $$
            This shows that $\Phi \colon Sz(T) \to T'$ is an inductive radiant.\\

            We observe that $\Phi \iota = \varphi$. Indeed, since $\varphi(x)^+ \leq \varphi(x^+)$ by \eqref{ip2}, and $\varphi(x^+) \in T'^+$ by Proposition \ref{prop:preradiant}, we have $\varphi(x^+) = \varphi(x^+)^+$ and $\varphi(x)^+ \wedge \varphi(x^+)^+ = \varphi(x)^+$. Therefore,
                $$ \Phi\iota(x) = \Phi(\{x,x^+\},x) = \left( \varphi(x)^+ \wedge \varphi(x^+)^+ \right) \varphi(x) = \varphi(x)^+\varphi(x) = \varphi(x), $$
          where the last identity follows from \eqref{c3}. 
          
          Moreover, $\Phi$ is the unique inductive radiant such that $\Phi \iota = \varphi$. Indeed, by Lemma \ref{lema:generated}, the Szendrei expansion $Sz(T)$ is generated by $\iota(T)$ under composition, range, and corestriction. Since inductive radiants preserve these operations, any two inductive radiants that agree on $\iota(T)$ must be equal. The converse follows from Lemma \ref{lema:composition}.
        \end{proof}
    \end{theorem}

\section{Second part of the ESN-type Theorem} In this section we prove that the category of left restriction semigroupoids and premorphisms is isomorphic to the category of li-constellations and inductive preradiants, thus establishing the second part of our ESN-type theorem.\\

    Throughout this section we let $S$ and $S'$ denote left restriction semigroupoids, and $T$ and $T'$ denote li-constellations. Recall that $C(S)$ denotes the li-constellation associated to $S$, as given in Lemma~\ref{lema:sgpd-to-const}, and that $S(T)$ denotes the left restriction semigroupoid associated to $T$, as constructed in Lemma~\ref{prop:const-to-sgpd}.

    \begin{lemma} \label{lema:premorphism1}
        If $\varphi \colon S \to S'$ is a premorphism, then $C\varphi = \varphi \colon C(S) \to C(S')$ is an inductive preradiant.

        \begin{proof}
            \eqref{ip1} If $x \bullet y$ is defined, then so is $xy$. By \eqref{pm1}, it follows that both $\varphi(x)\varphi(y)$ and $\varphi(x)^+ \varphi(xy)$ are defined and equal. From Proposition \ref{prop:objetos}, we then obtain:
                $$ \varphi(x) \otimes \varphi(y) = \varphi(x)\varphi(y) = \varphi(x)^+\varphi(xy) = \varphi(x)^+ \otimes \varphi(xy). $$
            In particular, this implies that $\varphi(x) | \varphi(y)^+ \neq \emptyset$ and $\varphi(x)^+ | \varphi(xy)^+ \neq \emptyset$.\\

            \noindent\eqref{ip2} This is \eqref{pm2}.\\

            \noindent\eqref{ip3} This corresponds to Lemma \ref{lema:pm}(b), since $S$ and $C(S)$ share the same partial order.\\

            \noindent\eqref{ip4} If $x | e \neq \emptyset$, then $xe$ is defined. From \eqref{lr4}, we have $xe = (xe)^+ x \leq x$. Therefore, by applying Lemma \ref{lema:pm}(b) and \eqref{wo2}, we obtain
                $$ xe \leq x \implies \varphi(xe) \leq \varphi(x) \implies \varphi(xe)^+ \leq \varphi(x)^+. $$
           On the other hand, since $xe$ is defined, it follows from \eqref{pm1} that $\varphi(x)\varphi(e)$ is also defined, and from Lemma \ref{lema:pm}(a) that $\varphi(e) \in S'^+$. Therefore, $\varphi(x) | \varphi(e) \neq \emptyset$, and we have:
            \begin{align*}
                \varphi(x)|\varphi(e) &= \varphi(x)\varphi(e) \\
                &= \varphi(x)^+ \varphi(xe) & \eqref{pm1} \\
                &= \varphi(x)^+ \varphi(xe)^+ \varphi(xe) & \eqref{lr1} \\
                &= \varphi(xe)^+ \varphi(xe) & (xe \leq x) \\
                &= \varphi(xe) & \eqref{lr1} \\
                &= \varphi(x|e).
            \end{align*}

            \noindent\eqref{ip5} If $e \leq x^+$, then $\varphi(e) \leq \varphi(x^+)$ by Lemma \ref{lema:pm}(b), and $\varphi(x)^+ \leq \varphi(x^+)$ by \eqref{pm2}. Therefore, $\varphi(e)$ and $\varphi(x)^+$ belong to the same connected component of $S'^+$, and by \eqref{wo9}, the corestriction $\varphi(e) | \varphi(x)^+$ is defined. Moreover, by Lemma \ref{lema:pm}(a) and Lemma \ref{lema:lr}(a), we have $\varphi(e) = \varphi(e)^+$. Hence,
            \begin{align*}
                \varphi(e) | \varphi(x)^+ &= \varphi(e)^+ \varphi(x)^+ \\
                &= (\varphi(e)^+\varphi(x))^+ & \eqref{lr3} \\
                &= (\varphi(e)\varphi(x))^+ \\
                &= (\varphi(e)^+ \varphi(ex))^+ & \eqref{pm1} \\
                &= \varphi(e)^+ \varphi(ex)^+ & \eqref{lr3} \\
                &= \varphi(e)\varphi(e|x)^+.
            \end{align*}
             On the other hand, we have:
             \begin{align*}
                 \varphi(e|x)^+ &= \varphi(ex)^+ \\
                 &\leq \varphi((ex)^+) & \eqref{pm2} \\
                 &= \varphi(ex^+) & \eqref{lr3} \\
                 &= \varphi(e). & (e \leq x^+)
             \end{align*}
             Therefore, $\varphi(e)\varphi(e|x)^+ = \varphi(e|x)^+$ and $\varphi(e)|\varphi(x)^+ = \varphi(e|x^+)$ as corestriction. This shows that $C\varphi = \varphi$ is an inductive preradiant from $C(S)$ to $C(S')$.
        \end{proof}
    \end{lemma}

    \begin{lemma} \label{lema:premorphism2}
        If $\varphi \colon T \to T'$ is an inductive preradiant, then $G\varphi = \varphi \colon G(T) \to G(T')$ is a premorphism.

        \begin{proof}
            By Theorem~\ref{teo:orderedSzendrei}, there exists a unique inductive radiant $\Phi \colon Sz(T) \to T'$ such that $\Phi \circ \iota = \varphi$, where $\iota \colon T \to Sz(T)$ is the inductive preradiant defined by $\iota(x) = {(x^+, x), x}$ for each $x \in T$.\\

            \noindent\eqref{pm1} If $s \otimes t$ is defined, then $s|t^+ \neq \emptyset$. By \eqref{ip4}, it follows that $\varphi(s)|\varphi(t^+) \neq \emptyset$. Moreover, since $\varphi(t)^+ \leq \varphi(t^+)$ by \eqref{ip2}, we deduce from \eqref{wo6} that $\varphi(s)|\varphi(t)^+ \neq \emptyset$. Hence, $\varphi(s) \otimes \varphi(t)$ is defined, and
            \begin{align*}
                \varphi(s) \otimes \varphi(t) &= \Phi\iota(s) \otimes \Phi\iota(t) \\
                &= (\Phi\iota(s)|\Phi\iota(t)^+)\Phi\iota(t) \\
                &= \Phi((\iota(s)|\iota(t)^+) \iota(t)) & \eqref{ir2}, \eqref{ir4}, \eqref{ir1} \\
                &= \Phi(\iota(s) \otimes \iota(t)) \\
                &= \Phi(\iota(s)^+ \otimes \iota(s \otimes t)) & \eqref{lema:iota} \\
                &= \Phi((\iota(s)^+|\iota(s \otimes t)^+)\iota(s \otimes t)) \\
                &= (\Phi\iota(s)^+|\Phi\iota(s \otimes t)^+) \Phi\iota(s \otimes t) & \eqref{ir1}, \eqref{ir4}, \eqref{ir2} \\
                &= \Phi\iota(s)^+ \otimes \Phi\iota(s \otimes  t) \\
                &= \varphi(s)^+ \otimes \varphi(s \otimes t).
            \end{align*}
            In particular, we obtain that $\varphi(s)^+ \otimes \varphi(s \otimes t)$ is defined. Note that \eqref{pm2} corresponds exactly to \eqref{ip2}. Therefore, $G\varphi = \varphi \colon G(T) \to G(T')$ is a premorphism.
        \end{proof}
    \end{lemma}

    \begin{prop}
        Li-constellations (as objects) and inductive preradiants (as morphisms) form a category, where the domain and range are given by the domain and range of functions, and composition is given by the usual composition of functions.

        \begin{proof}
            The proof is analogous to that of Proposition~\ref{prop:1st-category}, noting that Lemmas~\ref{lema:premorphism1} and~\ref{lema:premorphism2} establish a bijection between premorphisms and inductive preradiants, and that left restriction semigroupoids and premorphisms form a category by Proposition~\ref{prop:sgpd-cat-2}.
        \end{proof}
    \end{prop}

    We now establish the second part of our ESN-type theorem.

    \begin{theorem} \label{teo_princ_2}
        The category of left restriction semigroupoids and premorphisms is isomorphic to the category of li-constellations and inductive preradiants.

        \begin{proof}
         From Proposition~\ref{lema:sgpd-to-const} and Lemma~\ref{lema:premorphism1}, it follows that the functor $C$ maps the category of left restriction semigroupoids and premorphisms to the category of li-constellations and inductive preradiants. By Proposition~\ref{prop:const-to-sgpd} and Lemma~\ref{lema:premorphism2}, the functor $G$ maps the category of li-constellations and premorphisms to the category of left restriction semigroupoids and premorphisms. Moreover, $C$ and $G$ are mutually inverse, since
                $$ CG\varphi = C\varphi = \varphi \quad\text{and}\quad GC\varphi = G\varphi = \varphi, $$
            by definition, and $C(G(T)) = T$ and $G(C(S)) = S$ by Proposition \ref{prop:objetos}. It remains to verify that $C$ and $G$ are functors. Indeed, for every left restriction semigroupoid $S$, we have $Cid_S = id_S = id_{C(S)}$, and if $\varphi \colon S \to S'$ and $\psi \colon S' \to S''$ are morphisms, then
                $$ C(\psi \circ \varphi) = \psi \circ \varphi = C\psi \circ C\varphi, $$
            so $C$ is a functor. Analogously, for every li-constellation $T$, we have $Gid_T = id_T = id_{G(T)}$, and if $\varphi \colon T \to T'$ and $\psi \colon T' \to T''$ are inductive premorphisms, then $G(\psi \circ \varphi) = G\psi \circ G\varphi$. Therefore, $G = C^{-1}$, and the functors $C$ and $G$ establish an isomorphism between the respective categories.
        \end{proof}
    \end{theorem}

   Theorems~\ref{teo:1st-isomorphism} and~\ref{teo_princ_2} together yield the main result of this paper, stated below.

    \begin{theorem} \label{teo:principal}
        (ESN-type theorem for left restriction semigroupoids).
        \begin{itemize}
            \item[(i)] The category of left restriction semigroupoids and morphisms is isomorphic to the category of li-constellations and inductive radiants.
            \item[(ii)] The category of left restriction semigroupoids and premorphisms is isomorphic to the category of li-constellations and inductive preradiants.
        \end{itemize}
    \end{theorem}

\section{Particular Cases}

In this section we show that our one-sided ESN-type theorem unifies and generalizes results for three important classes of algebraic structures: left restriction categories, left restriction semigroups, and inverse semigroupoids. First we establish an ESN-type theorem for the class of left restriction categories. Then we show that Theorem~\ref{left-ESN}, which provides an ESN-type result for left restriction semigroups, follows as a particular case of our main result, Theorem~\ref{teo:principal}. Finally we derive a one-sided ESN-type theorem for the class of inverse semigroupoids, whose two-sided version is given in Theorem~\ref{multi-ESN}. We begin by examining additional properties on a li-constellation and identifying their corresponding formulations in the setting of left restriction semigroupoids.

    \begin{defi}
        A li-constellation $(T, \leq)$ is called:
\begin{itemize}
    \item[(ND)] \textit{non-degenerate} if, for every $x \in T$, there exists $e \in T^+$ such that $x|e \neq \emptyset$.

    \item[(LC)] \textit{locally complete} if each connected component $\omega(e)$ of $(T^+, \leq)$ has a maximum element, that is, an element $1 \in \omega(e)$ such that $x \leq 1$ for every $x \in \omega(e)$.

    \item[(U)] \textit{unitary} if it is locally complete and, for every $1 \in T^+$ that is the maximum element of some connected component of $(T^+, \leq)$, we have $x|1 = x$ whenever $x|1 \neq \emptyset$.
\end{itemize}
    \end{defi}

    We say that an element $e$ of a semigroupoid $S$ is a \textit{left identity} if $ee$ is defined and $es = s$ whenever $es$ is defined. Analogously, we say that $e$ is a \textit{right identity} if $ee$ is defined and $se = s$ whenever $se$ is defined. An \textit{identity} is an element that is both a left identity and a right identity.

Notice that, for a given $s \in S$, there is at most one identity $e$ such that $es$ is defined. Indeed, if $e, f \in S$ are identities and both $es$ and $fs$ are defined, then $es = e(fs)$ is defined, so $ef$ is defined. Since both $e$ and $f$ are identities, it follows that $e = ef = f$. Analogously, there is at most one identity $e$ such that $se$ is defined.

We recall that, in a left restriction semigroupoid, $ee$ is defined for every $e \in S^+$.

    \begin{lemma} \label{lema:particular}
        The li-constellation $C(S)$ associated to a left restriction semigroupoid is:
        \begin{itemize}
            \item[(a)] non-degenerate if and only if $S^s \neq \emptyset$ for every $s \in S$.
            \item[(b)] locally complete if and only if, for every $s \in S$, there exists a left identity $e \in S^+$ such that $es$ is defined.
            \item[(c)] unitary if and only if, for every $s \in S$, there exists an identity $e \in S$ such that $es$ is defined.
        \end{itemize}

        \begin{proof}
            (a) Let $s \in S$ and suppose that $t \in S^s$, that is, $st$ is defined. Then $st^+$ is defined, so $s|t^+ \neq \emptyset$. Conversely, suppose there exists $e \in C(S)^+$ such that $s|e \neq \emptyset$. Then $s \otimes e = (s|e) \bullet e$ is defined. Since the pseudo-product $\otimes$ coincides with the composition in $S$, it follows that $e \in S^s$. Therefore, $C(S)$ is non-degenerate if and only if, for every $s \in C(S)$, there exists $e \in C(S)^+$ such that $s|e \neq \emptyset$, which is equivalent to $S^s \neq \emptyset$ for every $s \in S$.\\

            (b) Let $s \in S$ and suppose that $C(S)$ is locally complete. Then there exists a maximum element $1 \in \omega(s^+)$.

            Since $1s^+$ is defined, it follows that $1s$ is also defined. We claim that $1 \in S^+$ is a left identity. Indeed, if $1t$ is defined, then $1t^+$ is also defined, so $t^+ \in \omega(1) = \omega(s^+)$. As $1$ is the maximum element of $\omega(s^+)$, we have $1t^+ = t^+$. Therefore, $1t = 1t^+t = t^+t = t$.

            Conversely, suppose that for every $s \in S$ there exists a left identity $1_s \in S^+$ such that $1_s s$ is defined. Let $\omega(e)$ be a connected component of $(C(S)^+, \leq)$. Since $1_e \in S^+$ and $1_e e$ is defined, we have $1_e \in \omega(e)$. We claim that $1_e$ is the maximum element of $\omega(e)$. Indeed, if $f \in \omega(e)$, then $1_e f$ is defined. As $1_e$ is a left identity, it follows that $1_e f = f$, hence $f \leq 1_e$.\\

            (c) Let $s \in S$ and suppose that $C(S)$ is unitary. By item (b), the maximum element $1 \in \omega(s^+)$ is a left identity and $1s$ is defined. We claim that $1$ is also a right identity. Indeed, if $t1$ is defined, then $t|1 \neq \emptyset$, and in this case we have $t1 = t \otimes 1 = (t|1) \bullet 1 = t|1 = t$.

            Conversely, suppose that for every $s \in S$ there exists an identity $1_s \in S$ such that $1_s s$ is defined. Since $1_s$ is, in particular, a right identity, it follows from \eqref{lr1} that $1_s = 1_s^+ 1_s = 1_s^+$, hence $1_s \in S^+$. By item (b), it follows that $C(S)$ is locally complete. Let $\omega(e)$ be a connected component of $(C(S)^+, \leq)$, and let $1_e$ be its maximum element. Then $1_e$ is an identity. Therefore, if $x \in C(S)$ is such that $x|1_e \neq \emptyset$, we have $x|1_e = x1_e = x$. This shows that $C(S)$ is unitary.
        \end{proof}
    \end{lemma}

   Now we present examples of left restriction semigroupoids whose associated li-constellations $C(S)$ satisfy the following: only condition (ND); none of the conditions; only conditions (LC) and (U); only conditions (ND) and (LC); and only condition (LC). Observe that, by definition, if condition (U) is satisfied, then condition (LC) is also satisfied.

    \begin{exe} \label{exe:6.3}
        Let $S = \{e, f, 0\}$ with $S^{(2)} = S \times S$, and define $xx = x$, $x^+ = x$ for every $x \in S$, and $ef = fe = 0e = e0 = f0 = 0f = 0$. Then $(S, +)$ is a left restriction semigroupoid. The li-constellation $C(S)$ satisfies condition (ND), since $S^x = S$ for every $x \in S$, but it does not satisfy condition (LC), as $S$ has no left identities.
    \end{exe}

    \begin{exe} \label{exe:6.4}
        Let $S = \{e,f,0\}$ be the left restriction semigroupoid with the structure given in Example \ref{exe:6.3}. Extend this structure to $S' = S \cup \{s\}$ by defining $es = fs = 0s = s$ and $s^+ = 0$. Then $(S', +)$ is a left restriction semigroupoid. The li-constellation $C(S')$ satisfies neither condition (LC) nor condition (ND), since $S^s = \emptyset$.
    \end{exe}

    \begin{exe} \label{exe:6.5}
        Let $S = \{x^+, x\}$ with $S^{(2)} = \{ (x^+, x^+), (x^+, x) \}$, and define $x^+ x^+ = x^+$ and $x^+ x = x$. Then $(S, +)$ is a left restriction semigroupoid. The li-constellation $C(S)$ satisfies condition (U), since $x^+$ is an identity and both $x^+ x^+$ and $x^+ x$ are defined. However, it does not satisfy condition (ND), as $S^x = \emptyset$.
    \end{exe}

    \begin{exe} \label{exe:6.6}
        Let $S = \{e,x,y,x^+,y^+\}$ be the semigroupoid with the multiplication table
        \begin{center}
            \begin{tabular}{c|c}
                 & $t$ \\\hline
                $s$ & $st$
            \end{tabular}\quad
            \begin{tabular}{c|c|c|c|c|c}
              & $e$ & $x$ & $y$ & $x^+$ & $y^+$ \\\hline
                $e$ & $e$ & - & - & - & - \\
                $x$ & $y$ & - & - & - & - \\
                $y$ & $y$ & - & - & - & - \\
                $x^+$ & - & $x$ & $y$ & $x^+$ & $y^+$ \\
                $y^+$ & - & $y$ & $y$ & $y^+$ & $y^+$ 
            \end{tabular}
        \end{center}
       Then $(S, +)$ is a left restriction semigroupoid with $S^+ = \{e, x^+, y^+\}$. The li-constellation $C(S)$ satisfies conditions (ND) and (LC), since $x^+$ and $e$ are left identities in $S^+$. However, it does not satisfy condition (U), as there is no identity $f \in S$ such that $fe$ is defined.
    \end{exe}

    \begin{exe} \label{exe:6.7}
        Let $S = \{e, x, y, x^+, y^+\}$ be the left restriction semigroupoid with the structure given in Example \ref{exe:6.6}. Extend this structure to $S' = S \cup \{s\}$ by defining $es = s$ and $s^+ = e$. Then $(S', +)$ is a left restriction semigroupoid. The li-constellation $C(S')$ satisfies condition (LC), since $e$ is a left identity for $s$. However, it does not satisfy conditions (ND) and (U), as $S'^s = \emptyset$. 
    \end{exe}

\subsection{Left restriction categories} In this subsection we characterize the non-degenerate and unitary (and therefore locally complete) li-constellations as precisely those arising from left restriction categories via the correspondence established in Theorem~\ref{teo:principal}. Consequently, we obtain an ESN-type theorem for left restriction categories.

    \begin{defi} \cite[Definition 2.1.1]{cockett2002restriction}
        Let $\mathcal{C} = (\mathcal{C}_0,\mathcal{C}_1,D,R,\circ)$ be a category. A \textit{left restriction} structure on $\mathcal{C}$ is a function $+ \colon \mathcal{C}_1 \to \mathcal{C}_1$ that assigns an arrow $f^+ \colon B \to B$ for each arrow $f \colon A \to B$ and such that the following conditions are satisfied:
        \begin{enumerate} \Not{R.}
            \item $f^+f = f$ for every $f$;
            \item $f^+ g^+ = g^+ f^+$ whenever $R(f) = R(g)$;
            \item $f^+ g^+ = (f^+ g)^+$ whenever $R(f) = R(g)$;
            \item $fg^+ = (fg)^+f$ whenever $D(f) = R(g)$.
        \end{enumerate}
        In this case, $(\mathcal{C},+)$ is called a \emph{left restriction category}.
    \end{defi}

    Every left restriction category is a left restriction semigroupoid. Therefore, to each left restriction category $(\mathcal{C}, +)$ we can associate a li-constellation $C(\mathcal{C})$. Recall that if $e \in \mathcal{C}_0$, then $e$ is an identity. In particular we have $e = e^+e = e^+$, so $e \in \mathcal{C}^+$.

    \begin{lemma} \label{lema:cat-const}
        If $(\mathcal{C},+)$ is a left restriction category, then $C(\mathcal{C})$ is a non-degenerate and unitary li-constellation. Conversely, if $(T,\leq)$ is a non-degenerate and unitary li-constellation, then $G(T)$ is a left restriction category.

        \begin{proof}
            If $(\mathcal{C}, +)$ is a left restriction category, then for every $s \in \mathcal{C}$, the composition $sD(s)$ is defined, and there exists an identity $R(s) \in \mathcal{C}$ such that $R(s)s$ is defined. Therefore, by Lemma~\ref{lema:particular}(a) and (c), it follows that $C(\mathcal{C})$ is non-degenerate and unitary.

            Conversely, suppose that $(T,\leq)$ is a non-degenerate and unitary (hence locally complete) li-constellation. Let $s \in G(T)$. Since $(T,\leq)$ is unitary and $C(G(T)) = T$, it follows from Lemma \ref{lema:particular}(c) that there is an identity $R(s) \in G(T)$ such that $R(s)s$ is defined. On the other hand, since $(T,\leq)$ is non-degenerate there is $e \in T^+$ such that $s|e \neq \emptyset$. As $(T,\leq)$ is locally complete, there exists a maximum element $D(s) \in \omega(e) \subseteq T^+$. From $e \leq D(s)$ and \eqref{wo6}, it follows that $s|D(s) \neq \emptyset$, so the composition $sD(s)$ is defined. Furthermore, $D(s)$ is an identity in $G(T)$, since for any $t \in G(T)$, the composition $tD(s)$ is defined if and only if $t|D(s) \neq \emptyset$, and in this case we have:
                $$ tD(s) = t \otimes D(s) = (t|D(s)) \bullet D(s) = t|D(s) = t, $$
            where the last identity follows from the fact that $(T,\leq)$ is unitary. Since there is at most one identity $R(s)$ such that $R(s)s$ is defined and at most one identity $D(s)$ such that $sD(s)$ is defined, we can define functions $D,R \colon G(T) \to G(T)$ that assign to each $s \in G(T)$ the unique identity $D(s)$ (respectively, $R(s)$) such that $sD(s)$ (respectively, $R(s)s$) is defined.

            The functions $D$ and $R$ endow $G(T)$ with a category structure. Indeed, observe that
            \begin{align*}
                (s,t) \in G(T)^{(2)} &\iff (sD(s),R(t)t) \in G(T)^{(2)} \\
                &\iff (D(s),R(t)) \in G(T)^{(2)} \iff D(s) = R(t).
            \end{align*}
            Moreover, if the composition $st$ is defined, then so are $s(tD(t)) = (st)D(t)$ and $(R(s)s)t = R(s)(st)$. It follows that $D(st) = D(t)$ and $R(st) = R(s)$. Hence, the quadruple $(G(T), D, R, \otimes)$ defines a category. Since $G(T)$ is already a left restriction semigroupoid, it follows that $G(T)$ is a left restriction category.
        \end{proof}
    \end{lemma}

 From the previous lemma we conclude that none of the left restriction semigroupoids presented in Examples~\ref{exe:6.3} through~\ref{exe:6.7} can be equipped with a category structure. 
 
 By combining Theorem~\ref{teo:principal} with Lemma~\ref{lema:cat-const} we obtain the following result:

    \begin{theorem} (ESN-type theorem for left restriction categories) \label{ESN-rCat}
        \begin{itemize}
            \item[(i)] The category of left restriction categories and restriction morphisms is isomorphic to the category of non-degenerate unitary li-constellations and inductive radiants.
            \item[(ii)] The category of left restriction categories and premorphisms is isomorphic to the category of non-degenerate unitary li-constellations and inductive preradiants.
        \end{itemize} 
    \end{theorem}

\subsection{Left restriction semigroups} In this subsection we explore the relationship between li-constellations and inductive left constellations, as defined in~\cite[Definition 2.1]{gould2009restriction}. We conclude that Theorem~\ref{teo:principal} generalizes the one-sided ESN-type theorem for left restriction semigroups developed by Gould and Hollings (Theorem~\ref{left-ESN}).\\

    A left restriction semigroup is a left restriction semigroupoid that is also a semigroup; that is, it satisfies $S^{(2)} = S \times S$. According to \cite[Theorem 4.13]{gould2009restriction} and \cite[Theorem 6.6]{gould2011actions}, left restriction semigroups are related via an ESN-like theorem (Theorem~\ref{left-ESN}) to inductive left constellations. By Proposition~\ref{prop:i-constellations}, the inductive left constellations are precisely the li-constellations $(T, \leq)$ for which $x|e \neq \emptyset$ for every $x \in T$ and $e \in T^+$.

    \begin{lemma} \label{lema:sg-const}
        The following conditions are equivalent.
        \begin{itemize}
            \item[(a)] $S$ is a left restriction semigroup.
            \item[(b)] $C(S)$ is non-degenerate, and $(C(S)^+,\leq)$ is a meet-semilattice.
            \item[(c)] For every $x \in C(S)$ and $e \in C(S)^+$, we have $x|e \neq \emptyset$.
        \end{itemize}

        \begin{proof}
            (a) $\implies$ (b). If $S$ is a semigroup, then $S^s = S \neq \emptyset$ for every $s \in S$. By Lemma~\ref{lema:particular}(a), it follows that $C(S)$ is non-degenerate. Moreover, since the product $ef$ is defined for all $e, f \in S^+$, we have $\omega(e) = S^+$ for every $e \in S^+$. Thus, $S^+$ consists of a single connected component, and consequently, $(C(S)^+, \leq)$ is a meet-semilattice. \\

            (b) $\implies$ (c). Let $x \in C(S)$ and $e \in C(S)^+$. Since $C(S)$ is non-degenerate, there exists $f \in C(S)^+$ such that $x|f \neq \emptyset$. As $(C(S)^+, \leq)$ is a meet-semilattice, the meet $e \wedge f$ exists in $C(S)^+$ and satisfies $e \wedge f \leq e$ and $e \wedge f \leq f$. Applying \eqref{wo6} twice, we conclude that $x|e \neq \emptyset$.\\

            (c) $\implies$ (b). Let $s, t \in S$. Since the pseudo-product $\otimes$ coincides with the composition in $S$, it follows that $st$ is defined if and only if $s \otimes t$ is defined, which in turn holds if and only if $s|t^+ \neq \emptyset$. By condition (c), this last condition is always satisfied. Therefore, $S^{(2)} = S \times S$, and $S$ is a semigroup.
        \end{proof}
    \end{lemma}

   Notice that the left restriction semigroupoids $(S',+)$ from Example~\ref{exe:6.4} and $(S,+)$ from Example~\ref{exe:6.5} are not semigroups, yet $S'^+ = S' \setminus \{s\}$ and $S^+ = \{x^+\}$ are meet-semilattices. In these cases, however, the associated constellations fail to satisfy condition (ND). 

   By combining Proposition~\ref{prop:i-constellations} with Lemma~\ref{lema:sg-const}(b) and (c), we conclude that the class of non-degenerate li-constellations $(T, \leq)$ for which $(T^+, \leq)$ is a meet-semilattice coincides exactly with the class of inductive left constellations described in Theorem~\ref{left-ESN}. Moreover, in this setting, the $(2,1)$-morphisms correspond exactly to morphisms. Consequently, the theorem below may be viewed as a reformulation of the one-sided ESN-type theorem for left restriction semigroups in the language of li-constellations.
   
   By combining Theorem~\ref{teo:principal} with Lemma~\ref{lema:sg-const}(a) and (b), we obtain the following result.

    \begin{theorem} (ESN-type theorem for left restriction semigroups)
        \begin{itemize}
            \item[(i)] The category of left restriction semigroups and morphisms is isomorphic to the category of non-degenerate li-constellations $(T,\leq)$ such that $(T^+,\leq)$ is a meet-semilattice and inductive radiants.
            \item[(ii)] The category of left restriction semigroups and premorphisms is isomorphic to the category of non-degenerate li-constellations $(T,\leq)$ such that $(T^+,\leq)$ is a meet-semilattice and inductive preradiants.
        \end{itemize} 
    \end{theorem}

\subsection{Inverse Semigroupoids} In this subsection we establish a one-sided ESN-type theorem for the class of inverse semigroupoids, whose two-sided version was developed by DeWolf and Pronk (see Theorem~\ref{multi-ESN}). In particular we generalize the correspondence between inductive left constellations and inverse semigroups, as presented in~\cite[Corollary 5.4]{gould2009restriction}.\\

    A \textit{regular semigroupoid} is a semigroupoid $S$ such that for every $s \in S$, there exists $t \in S$ with $(s, t), (t, s) \in S^{(2)}$, $sts = s$, and $tst = t$. In this case, the element $t$ is called a \textit{pseudo-inverse} of $s$. An \textit{inverse semigroupoid} is a regular semigroupoid in which every element has a unique pseudo-inverse, denoted by $s^{-1}$.

    \begin{obs} \label{obs:inv-sgpd}
        (1) Let $S$ be a semigroupoid such that for every $s \in S$, there exists $t \in S$ with $(s,t), (t,s) \in S^{(2)}$ and $sts = s$. Then $S$ is regular. Indeed, let $r = tst$. Since $(s, r), (r, s) \in S^{(2)}$, we compute:
            $$ srs = (sts)ts = sts = s \quad\text{and}\quad rsr = t(sts)tst = t(sts)t tst = r. $$
        That is, if $sts = s$, then $r = tst$ is a pseudo-inverse of $s$.

        (2) A regular semigroupoid is an inverse semigroupoid if and only if the set $E(S)$ of idempotent elements of $S$ is commutative in the following sense: if $e, f \in E(S)$ and $ef$ is defined, then $fe$ is also defined and $ef = fe$. A proof of this fact can be found in~\cite[Lemma 3.3.1]{liu2016free}.
    \end{obs}

    Given an inverse semigroupoid $S$, the function $s \mapsto s^+ = ss^{-1}$ endows $S$ with the structure of a left restriction semigroupoid. Consequently, to every inverse semigroupoid $S$, we can associate a li-constellation $C(S)$. A \textit{constellation with right inverses} is a left constellation $T = (T,T^{(2)},\star,+)$ such that for every $x \in T$, there exists an element  $x^{-1} \in T$ satisfying $xx^{-1} = x^+$.

    \begin{lemma} \label{lema:inv-sgpd}
        If $S$ is an inverse semigroupoid, then $C(S)$ is a li-constellation with right inverses. Conversely, if $(T,\leq)$ is a li-constellation with right inverses, then $G(T)$ is an inverse semigroupoid.

        \begin{proof}
            Suppose that $S$ is an inverse semigroupoid, and let $s \in S$. Observe that
                $$ (s^{-1})^+ = (s^{-1})(s^{-1})^{-1} = s^{-1}s. $$
            Hence, the composition $s(s^{-1})^+$ is defined and satisfies $s(s^{-1})^+ = s$. It follows that, in $C(S)$, we have $s \bullet s^{-1} = ss^{-1} = s^+$. This shows that $C(S)$ has right inverses.

            Conversely, suppose that $(T,\leq)$ is a li-constellation with right inverses. From Proposition \ref{prop:corestriction}, we know that if the product $xy$ is defined in $T$, then $x|y^+ = x$. In this case,
                $$ x \otimes y = (x|y^+)y = xy. $$
            Since $T$ has right inverses, for every $x \in G(T)$ there exists $x^{-1} \in G(T)$ such that $xx^{-1} = x^+$. It follows that $x^+x = (xx^{-1})x$ is defined and satisfies $(xx^{-1})x = x$. Since $G(T)$ is a semigroupoid and $(x \otimes x^{-1}) \otimes x$ is defined, it follows that both $x \otimes x^{-1}$ and $x^{-1} \otimes x$ are defined, and $x \otimes x^{-1} \otimes x = x$. By Remark~\ref{obs:inv-sgpd}(1), this shows that $G(T)$ is a regular semigroupoid. Observe that if $ax = bx$ in $T$, then $bx^+$ is defined, since both $bx$ and $xx^{-1}$ are defined. In this case,
                $$ a = ax^+ = a(xx^{-1}) = (ax)x^{-1} = (bx)x^{-1} = b(xx^{-1}) = bx^+ = b. $$
            That is, if $ax = bx$ in $T$, then $a = b$.

We now show that $E(G(T)) = G(T)^+$. Let $e \in E(G(T))$ be an idempotent element. Then $e|e^+ \neq \emptyset$ and $(e|e^+)e = e = e^+e$. From the previous argument, it follows that $e|e^+ = e^+ \in T^+$. By Lemma~\ref{lema:constellation}(a), we have $e|e^+ = (e|e^+)^+$, and applying Lemma~\ref{lema:corestriction}(a), we obtain:
                $$ e = (e|e^+)e = (e|e^+)^+e = e|e^+ = e^+. $$
            Thus, $e = e^+ \in G(T)^+$ for every $e \in E(G(T))$. Since $(T^+, \leq)$ is a meet-semilattice, it follows that $(E(G(T)), \otimes) = (G(T)^+, \otimes)$ is commutative. Then, since $G(T)$ is regular, it follows from Remark~\ref{obs:inv-sgpd}(2) that $G(T)$ is an inverse semigroupoid.
        \end{proof}
    \end{lemma}

By combining Theorem~\ref{teo:principal} with Lemma~\ref{lema:inv-sgpd} we obtain the following result:

    \begin{theorem} (One-sided ESN-type theorem for inverse semigroupoids) \label{ESN-inv-sgpd}
        \begin{itemize}
            \item[(i)] The category of inverse semigroupoids and morphisms is isomorphic to the category of li-constellations with right inverses and inductive radiants.
            \item[(ii)] The category of inverse semigroupoids and premorphisms is isomorphic to the category of li-constellations with right inverses and inductive preradiants.
        \end{itemize}
    \end{theorem}

  From Theorems~\ref{ESN-inv-sgpd} and~\ref{multi-ESN} we conclude that the category of li-constellations with right inverses and inductive radiants is isomorphic to the category of locally inductive groupoids and locally inductive functors.

   The correspondence between inverse semigroupoids and locally inductive groupoids established in Theorem~\ref{multi-ESN} restricts to a correspondence between inverse categories and locally complete inductive groupoids, that is, inductive groupoids satisfying condition (LC). The following result shows that, under the one-sided ESN-type theorem, inverse categories correspond precisely to those inverse semigroupoids $S$ whose associated li-constellation $C(S)$ satisfies condition (LC).

    \begin{prop}
        A left restriction semigroupoid $\mathcal{C}$ is an inverse category if and only if $C(\mathcal{C})$ is a locally complete li-constellation with right inverses.

        \begin{proof}
           From Theorems~\ref{ESN-rCat} and~\ref{ESN-inv-sgpd}, we obtain that $\mathcal{C}$ is an inverse category if and only if $C(\mathcal{C})$ is a non-degenerate and unitary li-constellation with right inverses. Every unitary li-constellation is locally complete. Therefore, it remains to prove that every locally complete li-constellation with right inverses is non-degenerate and unitary.

            Suppose that $T$ is a locally complete li-constellation with right inverses. Since $T$ has right inverses, $G(T)$ is an inverse semigroupoid. In particular, for every $s \in G(T)$, we have $s^{-1} \in G(T)^s$. By Lemma~\ref{lema:particular}(a), this implies that $G(T)$ is non-degenerate. Furthermore, the left restriction structure on $G(T)$ is given by $s^+ = ss^{-1}$. Using identity~\eqref{lr4}, we obtain that if $xy$ is defined in $G(T)$, then
                $$ xy^+ = (xy)^+x = (xy)(xy)^{-1}x = xyy^{-1}x^{-1}x = xy^+x^{-1}x. $$
            Since $T$ is locally complete, Lemma~\ref{lema:particular}(b) ensures that for every $s \in G(T)$, there exists a left identity $e \in G(T)^+$ such that $es$ is defined. Suppose that $xe$ is defined. Since $e = e^+$, it follows from the computation above, from the fact that $e$ is a left identity, and from identity~\eqref{lr1}, that
                $$ xe = xex^{-1}x = xx^{-1}x = x^+x = x.  $$
            Therefore, $e \in G(T)$ is an identity. By Lemma~\ref{lema:particular}(c), this implies that $T$ is unitary.
        \end{proof}
    \end{prop}
    
    \bibliographystyle{abbrvnat}
    \footnotesize{\bibliography{ref}}
\end{document}